\documentclass[11pt, reqno]{amsart}

\usepackage{amsmath,amssymb,amscd,amsthm,amsxtra}
\usepackage[dvips]{graphics,epsfig}

\allowdisplaybreaks[2]

\newtheorem{theorem}{Theorem} [section]

\newtheorem{lemma}[theorem]{Lemma}
\newtheorem{proposition}[theorem]{Proposition}
\newtheorem{remark}[theorem]{Remark}

\newtheorem{corollary}[theorem]{Corollary}

\DeclareMathOperator*{\supp}{supp}

\newcommand{\I}{\hspace{0.5mm}\text{I}\hspace{0.5mm}}
\newcommand{\II}{\text{I \hspace{-2.8mm} I} }
\newcommand{\III}{\text{I \hspace{-2.9mm} I \hspace{-2.9mm} I}}

\newcommand{\noi}{\noindent}
\newcommand{\Z}{\mathbb{Z}}
\newcommand{\R}{\mathbb{R}}
\newcommand{\C}{\mathbb{C}}
\newcommand{\T}{\mathbb{T}}

\newcommand{\al}{\alpha}
\newcommand{\dl}{\delta}
\newcommand{\Dl}{\Delta}
\newcommand{\eps}{\varepsilon}

\newcommand{\ld}{\lambda}

\newcommand{\s}{\sigma}
\newcommand{\ft}{\widehat}
\newcommand{\wt}{\widetilde}
\newcommand{\cj}{\overline}
\newcommand{\dx}{\partial_x}
\newcommand{\dt}{\partial_t}

\newcommand{\embeds}{\hookrightarrow}

\newcommand{\biggparen}[1]{\biggl(#1\biggr)}

\newcommand{\jb}[1]
{\langle #1 \rangle}

\makeatletter
\DeclareRobustCommand\invft[1]{{\mathpalette\@widecheck{#1}}}
\def\@widecheck#1#2{%
   \setbox\z@\hbox{\m@th$#1#2$}%
   \setbox\tw@\hbox{\m@th$#1%
      \widehat{%
         \vrule\@width\z@\@height\ht\z@
         \vrule\@height\z@\@width\wd\z@}$}%
   \dp\tw@-\ht\z@
   \@tempdima\ht\z@ \advance\@tempdima2\ht\tw@ \divide\@tempdima\thr@@
   \setbox\tw@\hbox{%
      \raise\@tempdima\hbox{\scalebox{1}[-1]{\lower\@tempdima\box\tw@}}}%
   {\ooalign{\box\tw@ \cr \box\z@}}}
\makeatother

\numberwithin{equation}{section}
\numberwithin{theorem}{section}

\begin{document}

%\baselineskip = 20pt
%\date{\today}

\title
[Brownian motions on modulation spaces] {Modulation spaces,  Wiener
amalgam spaces, and Brownian motions}

\author{\'Arp\'ad B\'enyi}
\author{Tadahiro Oh}

\address{\'Arp\'ad B\'enyi\\
Department of Mathematics\\
Western Washington University\\
516 High Street\\
Bellingham, WA 98226, USA}

\email{arpad.benyi@wwu.edu}

\address{Tadahiro Oh\\
Department of Mathematics\\
Princeton University\\
Fine Hall, Washington Rd\\
Princeton, NJ 08544-1000}

\email{hirooh@math.princeton.edu}



\begin{abstract}
We study the local-in-time regularity of the Brownian motion with
respect to localized variants of modulation spaces $M^{p, q}_s$ and Wiener amalgam
spaces $W^{p, q}_s$. We show that the periodic Brownian motion belongs
locally in time to $M^{p, q}_s (\T)$ and $W^{p, q}_s (\T)$ for $(s-1)q < -1$,
and the condition on the indices is optimal. Moreover, with the
Wiener measure $\mu$ on $\T$, we show that $(M^{p, q}_s (\T), \mu)$ and
$(W^{p, q}_s (\T), \mu)$ form abstract Wiener spaces for the same range
of indices, yielding large deviation estimates. We also establish
the endpoint regularity of the periodic Brownian motion with respect to a
Besov-type space $\ft{b}^s_{p, \infty} (\T)$. Specifically, we prove that
the Brownian motion belongs to $\ft{b}^s_{p, \infty} (\T)$ for $(s-1) p =
-1$, and it obeys a large deviation estimate. Finally, we revisit
the regularity of Brownian motion on usual local Besov spaces $B_{p,
q}^s$, and indicate the endpoint large deviation estimates.

\end{abstract}

\maketitle

\tableofcontents

\section{Introduction}

Modulation spaces were born during the early eighties in pioneering
work of H. Feichtinger. In subsequent fruitful collaborations with
K. Gr\"ochenig \cite{FG1}, \cite{FG2} they established the basic
theory of these function spaces, in particular their invariance,
continuity, embeddings, and convolution properties. In contrast with
the Besov spaces, which are defined by a dyadic decomposition of the
frequency space, modulation spaces (and Wiener amalgam spaces) arise
from a uniform partition of the frequency space. Their appeal is due
to the fact that they can effectively capture the time-frequency
concentration of a distribution. Both modulation and Wiener amalgam
spaces are modeled on Lebesgue or Fourier-Lebesgue spaces, but they
are more flexible in that they allow a separate control of the local
regularity and the decay at infinity of a function. The central idea
in this time-frequency analysis is to start with a function that is
smoothly localized by a so-called window function and then take its
Fourier transform. The resulting short-time Fourier transform (STFT)
is also known, depending on its specific application and up to some
normalization, as the Gabor transform, ambiguity function, coherent
state transform, or Wigner distribution. Modulation spaces thus
capture the joint time-frequency concentration of a function by
appropriate decay and integrability conditions on the STFT. They
originate in the early foundations of quantum mechanics and
information theory. Engineers and physicists have sensed a huge
potential here before mathematicians did, and the past years are
seeing a resurgence in the applications of time-frequency analysis
to signal analysis, image processing, or information theory. Their rapid
development was almost in sync with that of wavelet theory.
Grossmann and Morlet \cite{grosmor}, for example, have realized
early on that both the wavelet transform and the STFT are special
cases of square integrable representations. Unsurprisingly perhaps,
since the STFT is a joint representation of a function in both time
and frequency, uncertainty principles can be thought of as
appropriate embeddings of modulation spaces; see Galperin and
Gr\"ochenig's article \cite{galpgro}. Returning to the original
analogy with the Besov spaces and their wavelet bases, modulation
spaces do not admit orthonormal bases of time-frequency shifts, but
rather so-called frames of time-frequency shifts. Due to Daubechies'
insightful work \cite{daub}, the concept of frame has quickly become crucial
in signal analysis. Within pure and applied mathematics, modulation
spaces and Wiener amalgam spaces are nowadays present in
investigations that concern problems in numerical analysis, operator
algebras, localization operators, Fourier multipliers,
pseudodifferential operators, Fourier integral operators,
non-linear partial differential equations, and so on. The list is
extensive, and we cannot hope to acknowledge here all those who made
the theory of modulation spaces such a successful story. Simply to
give a flavor of some of the recent works in partial differential
equations employing these spaces, we mention the contribution of
B\'enyi-G\"ochenig-Okoudjou-Rogers \cite{bgor} on unimodular
multipliers and the phase-space concentration of the solutions to
the free Schr\"odinger and wave equations, the work on the well-posedness of
non-linear (Sch\"odinger, Ginzburg-Landau, Klein-Gordon, KdV) equations with rough
data by B\'enyi-Okoudjou \cite{beko} and Wang et all \cite{balibo}, \cite{bahu},
\cite{bahud}, \cite{balichu}, or the articles of Cordero-Nicola on Strichartz estimates and the Schr\"odinger
equation with quadratic Hamiltonian \cite{cornic1}, \cite{cornic2}.

\newpage
This article lies at the interface between time-frequency analysis,
probability theory, and partial differential equations (PDEs). Its
main aim is to serve as a first bridge between these areas by
specifically pointing out the role that the spaces of time-frequency
analysis (modulation and Wiener amalgam) can play here. A central
question in applied mathematics and theoretical physics is how
initial data are propagated by non-linear PDEs. 
One important question is existence of an invariant measure for the flow.
Following Lebowitz, Rose, and Speer \cite{LRS},
Bourgain \cite{BO4, BO5} constructed invariant Gibbs measures
for some Hamiltonian PDEs on $\T$.
Such Gibbs measures can be regarded as weighted Wiener measures on $\T$ (= Brownian motion on $\T$), 
which are supported on function spaces with low regularity.
Therefore, in constructing a flow on the support of the Gibbs measure, 
one is often forced
to go beyond the usual Sobolev spaces, and cross over in the
more exotic realm of the variants of Fourier-Lebesgue spaces
discussed above. 
For example, in the case of the derivative non-linear
Schr\"odinger equation (DNLS), the flow is not well defined in the classical
$H^s(\T)$ space for $s<1/2$, yet almost sure global well-posedness can be established
by employing Fourier-Lebesgue spaces $\mathcal{F}L^{s, q}(\T)$
for some $s$ and $q$ with $(s-1)q < -1$; see the recent work of Nahmod, Oh,
Rey-Bellet, and Staffilani on DNLS \cite{nors}. 
See also the
earlier works of Bourgain on mKdV \cite{BO4}, Zakharov system \cite{BO5}
and Oh \cite{OH3, OHSBO}
for the use of Fourier-Lebesgue spaces $\mathcal{F}L^{s, \infty}(\T)$.

%
%Since in virtually
%all physical problems smooth data is perturbed by noise, in
%practice, the initial data hahave low regularity. A manifestation of
%equilibrium states for such PDEs is the existence of invariant
%measures under the PDE flow. A notion of closeness on probability
%measures of functions is needed when a PDE exhibits more than one
%invariant measure, and this notion of closeness is determined by
%regularity. 
%
%There is an increasing interest in the study of PDEs
%with random data in low regularity setting. 

We will study the Brownian motion, which is arguably the most
accessible continuous-time stochastic process yet it plays a central
role in both pure and applied mathematics. Its applications range
from the study of continuous-time martingales,  stochastic calculus
to control theory and financial mathematics. It is well known that
the integral of a Gaussian white noise is represented through
Brownian motion, and white noise is supported on functions of low
regularity. Thus, understanding the regularity of Brownian motion is
a natural problem. In the context of the classical function spaces
of PDEs (Sobolev and Besov), the local regularity of the Brownian
motion is well understood; see the works of Ciesielski \cite{CSki1},
\cite{CSki2} and Roynette \cite{ROY}. These results are summarized
in Subsection 2.1. An extension to the vector valued setting can be
found in Hynt\"onen and Veraar's paper \cite{hytver}. Some immediate
applications of Brownian motion regularity to stochastic integrals
and equations have appeared early on in the work of Ciesielski,
Kerkyacharian, and Roynette \cite{CSki3}. Schilling \cite{schil} has
extended the local and global regularity properties in \cite{CSki3},
as well as the ones for sample paths of Feller processes.
Interestingly, \cite{schil} presents a nice connection to the
asymptotic behavior of the symbol of the pseudodifferential operator
given in terms of the infinitesimal generator of the Feller process.
For a survey of recent advances in stochastic calculus with respect
to (fractional) Brownian motion and its connection to Malliavin
calculus, see Nualart's article \cite{nua}. In this work, we
establish the local regularity of the Brownian motion on appropriate
modulation spaces and Wiener amalgam spaces, and prove that it obeys
so called large deviation estimates. For the ``end-point'' results
we appeal to a Besov-type space introduced by Oh in a series of
works that were concerned with the invariance of the white noise for
the KdV equation \cite{OH4} and the stochastic KdV equation with
additive space-time white noise \cite{OH6}. Incidentally, our
results also recover the regularity of the white noise by dropping
one regularity from the one of the Brownian motion. We also revisit
the local regularity of Brownian motion on the usual Besov spaces. A
common thread throughout this work is the use of random Fourier
series.

\medskip

\noindent
{\bf Acknowledgments:}
The authors would like to thank the anonymous referees for thoughtful comments that have significantly improved this paper.
The second author would also like to thank Professor Jeremy Quastel for helpful discussions.

\subsection{Function spaces of time-frequency analysis}

We start by recalling some basic definitions regarding the function
spaces used throughout this work; see, for example, Gr\"ochenig's
book \cite{Gr}. Given a (fixed, non-zero) window function
$\phi\in\mathcal S(\R^d)$, the short-time Fourier transform (STFT)
$V_\phi f$ of a tempered distribution $u$ is
\[V_\phi u(x, \xi)=\langle u, M_\xi T_x \phi\rangle=\mathcal
F(u\overline{T_x\phi})(\xi) = \int_{\R^d} e^{-i y\cdot \xi}
\cj{\phi(y - x)} u (y) dy.\]

\noi Here, $\mathcal F(g)(\xi)=\widehat
g(\xi)=\displaystyle\int_{\mathbb R^d} g(x) e^{-ix\cdot \xi} dx$
denotes the Fourier transform of a distribution $g$,
while, for $x, \xi\in \R^d$, $M_\xi g(x)=e^{ix\cdot \xi}g(x)$, $T_x g(y)=g(y-x)$
denote the modulation and translation operators, respectively. $\mathcal F^{-1}(g)=\invft{g}$ will denote the inverse Fourier transform of $g$.
For $s\in\R$, we let $\langle \xi\rangle^s=(1+|\xi|^2)^{s/2}$.

The (continuous weighted) {\it modulation space} $M^{p,q}_s(\R^d),
1\leq p,q\leq\infty$,  consists of all tempered distributions
$u\in\mathcal S(\R^d)$ such that
$$V_\phi u (x, \xi)\langle \xi\rangle^s\in L^p_x L^q_\xi,$$

\noi and we equip the space $M^{p, q}_s(\R^d)$ with the norm
$$\|u\|_{M^{p,q}_s(\R^d)}=\big\|\|V_\phi u (x, \xi)\langle
\xi\rangle^s\|_{L^p_x(\R^d)}\big\|_{L^q_\xi(\R^d)}=
\biggparen{\int_{\R^{d}}\langle \xi\rangle^{sq}\biggparen{\int_{\R^{d}}|V_{\phi}u (x,
\xi)|^{p}\, dx}^{q/p}\, d\xi}^{1/q},$$
with the obvious modifications if $p=\infty$ or $q=\infty$. It is not hard to see that these spaces are Banach spaces and
two different windows yield equivalent norms, Moreover, the
classical Bessel potential spaces (or the usual $L^2$-based Sobolev spaces) are particular cases of modulation
spaces, $M^{2, 2}_s(\R^d)=H^s(\R^d)$, the duality and embeddings are natural
$(M^{p,q}_s(\R^d))'=M^{p',q'}_{-s}(\R^d)$, $M^{p_1,
q_1}_s(\R^d)\subset M^{p_2, q_2}_s(\R^d)$ for $p_1\leq p_2$ and
$q_1\leq q_2$, while $\mathcal S(\R^d)$ is dense in $M^{p,
q}_s(\R^d)$. The complex interpolation of these spaces is equally natural:
$(M^{p_1, q_1}_{s_1}, M^{p_2, q_2}_{s_2})_{[\theta]}=M^{p,q}_s$, where
$1/p=(1-\theta)/p_1+\theta/p_2, 1/q=(1-\theta)/q_1+\theta/q_2$ and $s=(1-\theta)s_1+\theta s_2$.

Closely related to these spaces are the (continuous weighted) {\it
Wiener amalgam spaces} $W^{p,q}_s(\R^d)$ which are now equipped with
the norm
$$\|u\|_{W^{p,q}_s(\R^d)}=\big\|\|V_\phi u (x, \xi)\langle
\xi\rangle^s\|_{L^q_\xi(\R^d)}\big\|_{L^p_x(\R^d)}=\biggparen{\int_{\R^{d}}\biggparen{\int_{\R^{d}}\langle \xi\rangle^{sq}|V_{\phi}u (x,
\xi)|^{q}\, d\xi}^{p/q}\, dx}^{1/p}.$$

\noi It is clear then that, in fact, $M^{p, q}_0(\R^d)=\mathcal F
W^{q, p}_0(\R^d)$. As such, completely analogous comments to the one
following the definition of modulation spaces hold for Wiener
amalgam spaces as well. For example, $W^{p_1, q_1}_s(\R^d)\subset
W^{p_2, q_2}_s(\R^d)$ for $p_1\leq p_2$ and $q_1\leq q_2$, $\mathcal
S(\R^d)$ is dense in $W^{p, q}_s(\R^d)$, and so on. Moreover, $M^{p,
q}_s(\R^d)\subset W^{p,q}_s(\R^d)$ for $q\leq p$, while the reverse
inclusion holds if $p\leq q$.

The (weighted) {\it Fourier-Lebesgue spaces} $\mathcal F L^{s,
p}(\R^d)$ are defined via the norm
\begin{equation}\label{WAmal1}
\|u\|_{\mathcal{F}L^{s, p}(\R^d)} = \|\jb{\xi}^s \ft{u}(\xi)
\|_{L_\xi^p(\R^d)}=\biggparen{\int_{\R^{d}}\langle \xi\rangle^{sp}|\ft{u}(\xi)|^p\, d\xi}^{1/p}.
\end{equation}

For our purposes, we will use the following equivalent definitions
for the norms of modulation and Wiener amalgam spaces. Let
$\psi \in \mathcal{S}(\R^d)$ such that $\supp \psi \subset [-1,
1]^d$ and $\sum_{n \in \Z^d} \psi(\xi -n) \equiv 1$. Then,
\begin{equation} \label{mod1}
\|u\|_{M^{p, q}_s(\R^d)} = \big\| \jb{n}^s \|\psi(D-n) u
\|_{L_x^p(\R^d)} \big\|_{\ell^q_n(\mathbb{Z}^d)}
\end{equation}
and
\begin{equation} \label{WAmal2}
\|u\|_{W^{p, q}_s(\R^d) } = \big\|  \| \jb{n}^s \psi(D-n) u
\|_{\ell^q_n(\mathbb{Z}^d) } \big\|_{L_x^p(\R^d)}.
\end{equation}
Here, we denoted
$$\psi (D-n) u(x)=\mathcal F^{-1}(\psi (\cdot-n)\mathcal F u(\cdot))(x)=
(2\pi)^{-d}\int_{\R^d} \psi (\xi-n)\widehat u (\xi)e^{i\cdot \xi}\, d\xi.$$
Thus, the modulation norm in (\ref{mod1}) can be spelled out, up to a $(2\pi)^{-d}$ factor, as
$$\biggparen{\sum_{n\in\mathbb{Z}^d}\biggparen{\int_{\R^d}\Big|\int_{\R^d}\psi (\xi-n)\ft{u}(\xi)e^{ix\cdot\xi}\, d\xi\Big|^p\,
dx}^{q/p}}^{1/q}.$$
One should contrast these definitions with the one of {\it Besov
spaces}. Let $\varphi_0, \varphi \in \mathcal{S}(\R^d)$ such that
$\supp \varphi_0 \subset \{ |\xi| \leq 2\}$, $\supp \varphi \subset
\{ \frac{1}{2}\leq |\xi| \leq 2\}$, and $ \varphi_0(\xi) + \sum_{j =
1}^\infty \varphi(2^{-j}\xi) \equiv 1.$ With $\varphi_j(\xi) =
\varphi(2^{-j}\xi),$ $ j\geq 1,$ we define the usual Besov space $B^s_{p, q}$
via the norm
\begin{equation} \label{besov1}
\|u\|_{B^s_{p, q}(\R^d) } = \big\| 2^{jsq} \|\varphi_j(D) u
\|_{L^p(\R^d)} \big\|_{\ell^q_j(\mathbb N)}
=\biggparen{\sum_{j\geq 0}2^{sjq}\|u*\invft{\varphi}_j\|_{L^p (\R^d)}^q }^{1/q}.
\end{equation}
Clearly, when $q=\infty$, the norm is modified into $\sup_{j\geq 0} 2^{sj}\|u*\invft{\varphi}_j\|_{L^p (\R^d)}.$
The embeddings between Besov, Sobolev and modulation spaces are well understood. For example, with some explicitly found
indices $s_1 (n, p, q), s_2 (n, p, q)$, we have $B^s_{p, q}(\R^d)\embeds M^{p, q}(\R^d)$ if and only if $s\geq s_1$ and
$M^{p, q}(\R^d)\embeds B^s_{p,q}(\R^d)$ if and only if $s\leq s_2$; see Sugimoto-Tomita \cite{suto}, also Okoudjou \cite{Ok}, Toft \cite{To}, and
Kobayashi-Sugimoto \cite{kosu}.

As stated in the introduction, our goal is to investigate the local-in-time regularity of the
Brownian motion. Without loss of generality, we can restrict ourselves to study the periodic Brownian motion;
see Subsection 2.2. In Subsection 2.1, we first summarize some of the known regularity results for Sobolev and Besov spaces.
Then, in Subsection 2.2, we establish similar results for modulation and Wiener amalgam spaces on the torus. We also observe that the same results
can be obtained on appropriately defined localized spaces, and that straightforward modifications allow us to
extend them for the Brownian motion on the $d$-dimensional torus. Moreover, in Section 3,
we will do much more by showing a nice connection to abstract Wiener spaces, and, in particular, their large deviation
estimates. In Subsection 2.3, we reprove some known results
for Besov spaces by working, as for the other spaces of
time-frequency analysis, on the Fourier side (compared to previous
proofs which were done on the physical side) and using a multilinear analysis argument. From this perspective,
the use of Fourier-Wiener series (see \eqref{2FW} below) can be
viewed as a unifying theme of this article. Since we are only
concerned with local regularity of (periodic) Brownian motions, we need a
``localized'' version of the time-frequency spaces defined above. These localizations
were first introduced and studied by Ruzhansky, Sugimoto,
Toft, and Tomita in \cite{RSTT}. Naturally, the spaces can be localized either in space (by
multiplying an object in the functional space by a smooth cut-off function) or in frequency (by multiplying the
Fourier transformation of an element in the functional space by a smooth cut-off). Proposition 2.1
and Remark 4.2 in \cite{RSTT} show that the localized versions, say in space, of the $M^{p,q}_s$, $W^{p, q}_s$ and $\mathcal FL^{q, s}$ spaces
coincide, with equivalence of norms. In particular, the modulation, Wiener amalgam, and Fourier-Lebesgue spaces on the torus are all the same.
Specifically, if for some $\psi$ with compact support in the discrete topology of $\Z^d$, we let
\begin{equation}\label{modtorus}
\|u\|_{M^{p, q}_s(\T^d)} = \big\| \jb{n}^s \|\psi(D-n) u
\|_{L_x^p(\T^d)} \big\|_{\ell^q_n(\mathbb{Z}^d)}
\end{equation}
and
\begin{equation} \label{WAmaltorus}
\|u\|_{W^{p, q}_s(\T^d) } = \big\|  \| \jb{n}^s \psi(D-n) u
\|_{\ell^q_n(\mathbb{Z}^d) } \big\|_{L_x^p(\T^d)},
\end{equation}
then, for all $1\leq p, q\leq \infty$, we have
$$M^{p, q}_s(\T^d)=W^{p, q}_s(\T^d)=\mathcal F L^{s, q}(\T^d)=\mathcal F \ell^{q}_s.$$

We will also use the Fourier-Besov spaces $\widehat{b}^s_{p, q}(\T^d)$. With the notation $\varphi_j$ introduced above for the
classical Besov spaces, the Fourier-Besov space norm is defined by:
\begin{equation} \label{Besov2}
 \|u\|_{\ft{b}^s_{p, q}(\T^d)} = \big\| \|\jb{n}^s \varphi_j(n)\ft{u}(n)\|_{L^p (\mathbb Z^d)} \big\|_{\ell^q_j(\mathbb N)}=
 \biggparen{\sum_{j\geq 0}\biggparen{\sum_{n\in\mathbb Z^d}\jb{n}^{sp}|\varphi_j (n)\widehat u(n)|^p}^{q/p}}^{1/q}.
\end{equation}
Clearly, this norm is equivalent to
\begin{equation}\label{Besov3}
\biggparen{\sum_{j\geq 0}\biggparen{\sum_{|n|\sim 2^j}2^{jsp}|\widehat u(n)|^p}^{q/p}}^{1/q}.
\end{equation}

\noi In particular, in our end-point regularity analysis we will
need the space $\ft{b}^s_{p, \infty}(\T^d)$ defined via
\begin{equation} \label{Besov}
\|u\|_{\widehat{b}^s_{p, \infty}(\T^d)}=
\sup_{j\geq 0}\|\jb{n}^s\varphi_j (n)\widehat{u}(n)\|_{L^p (\mathbb Z^d)} \sim\sup_{j\geq 0}
\bigg(\sum_{|n|\sim 2^j}2^{jsp}|\widehat{u}(n)|^p\bigg)^{1/p},
\end{equation}
where, for a given $n\in \mathbb Z^d$, $|n|\sim 2^j$ denoted $2^{j-1}<|n|\leq 2^j.$ This
space is well-suited for the analysis of the white noise, since, in particular, for $sp\leq-d$,  it contains its full support.
For $d=1$, this analysis of the white noise was initiated by Oh in \cite{OH4}, while the natural extension to any dimension $d$ can be
found in Veraar's work \cite{ver1}. 
Noting that the white noise is a derivative of the Brownian motion, 
some of the results in \cite{ver1} can be easily deduced from the
ones presented here; see Subsection 2.2. 
Lastly, note that Hausdorff-Young inequality immediately implies that $B^{s}_{p',
\infty}(\T^d)\subset \widehat{b}^s_{p, \infty}(\T^d)$ for $p>2$, and that the
equality holds for $p=2$, where $B^s_{p', \infty}(\T^d)$ is the usual periodic Besov space and $p'=p/(p-1)$ denotes 
the H\"older conjugate exponent of $p$.
Note also that $\mathcal F L^{s, p}(\T^d)=\widehat{b}^s_{p, p}(\T^d)$ for $p\geq 1$.

\subsection{Brownian motion}

Let $T>0$ and $(\Omega, \mathcal F, Pr)$ be a probability space. We define the
{\it $\R^n$-valued Brownian motion} (or Wiener process) as a
stochastic process $\beta: [0, T]\times \Omega\rightarrow \R^n$ that
satisfies the following properties:
\begin{itemize}
\item[(i)] $\beta(0) = 0$ almost surely (a.s. $\omega\in \Omega$)
\item[(ii)] $\beta(t)$ has independent increments, and
$\beta(t) - \beta(t')$ has the normal distribution with mean 0 and
variance $t-t'$ (for $0 \leq t' \leq t$).
\end{itemize}
Here, we abused the notation and wrote, for a given $t\in [0, T]$,
$\beta(t): \Omega\rightarrow \R^n$, $\beta(t)(\omega)=\beta (t,
\omega)$. This further induces a mapping $\Phi_\beta$ from $\Omega$
into the collection of functions from $[0, T]$ to $\R^n$ by
$(\Phi_\beta (\omega))(t)=\beta(t)(\omega)$. The law of the Brownian
motion, that is, the pushforward measure $(\Phi_\beta)_{*}(Pr)$ on
Borel sets in $(\R^n)^{[0, T]}$, is nothing but the classical
Wiener measure. Replacing ${\mathbb R}^n$ with some other Banach
space will take us into the realm of abstract Wiener spaces. It is
worth noting that, while there are several technicalities involved
in properly defining a Wiener measure in this general context, the
same point of view as in the classical $\mathbb R^n$ case is
employed.

\section{Regularity of Brownian motion}

In what follows, we first discuss some known results about the
regularity of Brownian motion on function spaces that are often used
in PDEs, that is, Sobolev and Besov spaces. All spaces are
considered local in time.

\subsection{Modulus of continuity}
\noi It is well-known that $\beta(t)$ is almost surely (a.s.)
continuous. Moreover, we have L\'evy's theorem for the modulus of
continuity:
\begin{equation}\label{LEVY}
 \limsup_{\substack{t - t'\leq \eps \downarrow 0\\ 0\leq t'< t\leq 1}}
\frac{|\beta(t) -\beta(t')|}{\sqrt{-2\eps \log \eps}} = 1,  \,
\mbox{a.s.}
\end{equation}

\noi It follows from \eqref{LEVY} that the Brownian motion is
a.s.~locally H\"older continuous of order $s$ for every $s
< 1/2$, and is a.s.~nowhere locally H\"older continuous of order $s$
for $s \geq 1/2$; see Revuz and Yor's book \cite{RY}.

\medskip
\noi {\it Regularity on Sobolev spaces $H^s (0, 1)$:}
By definition of the Brownian motion, we have
\[ \mathbb{E} \big[|\beta(t) - \beta(t')|^2\big]
 = \mathbb{E} \big[|\beta(t - t') |^2\big] = |t-t'|,\]
where $\mathbb{E}$ denotes expectation.
Thus, we have
\begin{align*}
\mathbb{E} \big[ \|\beta(t) \|_{\dot{H}^s(0, 1)}^2 \big] & =
\mathbb{E} \bigg[\int_0^1 \int_0^1 \frac{|\beta(t) -
\beta(t')|^2}{|t-t'|^{1+2s}} d t' dt \bigg]
= \int_0^1 \int_0^1 \frac{1}{|t-t'|^{2s}} d t' dt\\
& = 2 \int_0^1 \int_0^t \frac{1}{(t-t')^{2s}} d t' dt < \infty,
\end{align*}

\noi if and only if $ s < \frac{1}{2}$. Hence, the Brownian motion
$\beta(t)$ belongs almost surely to $H^s_\text{loc}$ for $ s <
\frac{1}{2}$. Indeed, one can also show that $\beta(t) \notin
H^s_{\text{loc}}$ a.s.~for $s \geq \frac{1}{2}$.

\medskip
\noi {\it Regularity on Sobolev spaces $\mathcal W^{s, p}(0, 1)$:}

Recall that we have, for $p> 0$,
\[ \mathbb{E} \big[ |\beta(t) - \beta(t')|^p \big]
=C_p |t - t'|^\frac{p}{2},
\]

\noi since $\beta(t) - \beta(t')$ is a mean-zero Gaussian with
variance $|t - t'|$. Then, by the characterization of $\dot{\mathcal
W}^{s, p}$ via the $L^p$ modulus of continuity  as in Tartar's book
\cite{Tartar}, we have
\begin{align*}
\mathbb{E} \big[ \|\beta(t) \|_{\dot{\mathcal W}^{s, p}(0, 1)}^p
\big] & = \mathbb{E} \bigg[\int_0^1 \int_0^1 \frac{|\beta(t) -
\beta(t')|^p}{|t-t'|^{1+sp}} d t' dt \bigg]
\sim \int_0^1 \int_0^1 |t-t'|^{-1- sp+ \frac{p}{2}} d t' dt\\
& = 2 \int_0^1 \int_0^t |t-t'|^{-1- sp+ \frac{p}{2}} d t' dt <
\infty,
\end{align*}

\noi if and only if $ s < \frac{1}{2}$. Hence, the Brownian motion
$\beta(t)$ belongs almost surely to $\mathcal W^{s, p}_\text{loc}$
for $ s < \frac{1}{2}$.

\medskip
\noi {\it Regularity on Besov spaces $B^s_{p, q}(0, 1)$:}

Ciesielski \cite{CSki1}, \cite{CSki2} and Roynette \cite{ROY} proved
that if $s < \frac{1}{2}$, then the Brownian motion $\beta (t)$
belongs a.s. to $B^s_{p, q}(0, 1)$ for all $p, q \geq 1$, and that
if $s > \frac{1}{2}$, then $\beta (t)\notin B^s_{p, q}(0, 1)$ a.s.
for any $p, q \geq 1$. Regarding the endpoint regularity of the
Brownian motion, it was also shown that for $s = \frac{1}{2}$, the
Brownian motion $\beta (t) \in B^\frac{1}{2}_{p, q}(0, 1)$ a.s. if and only if
$1 \leq p < \infty$ and $q = \infty$. Moreover, if the latter holds,
there exists $c_p>0$ such that $\| \beta (t) \|_{B^\frac{1}{2}_{p,
\infty}(0, 1)} \geq c_p$ a.s. The proof is based on the Schauder
basis representation (or Franklin-Wiener series) of the Brownian
motion; see Kahane's book \cite{Kahane}.

In Subsection 2.3, we  will present an alternate  proof of these
results which is of interest in its own, by  using random Fourier
series (also referred to as Fourier-Wiener series, see \eqref{2FW} below.) It is worth noting that
the independence on the index $p<\infty$ of the regularity result on, say,
$\mathcal W^{s, p}_\text{loc}$ is unsurprising from the perspective of random Fourier series.
Without delving into technicalities, this follows essentially by observing that,
if $F(t)=\sum_{n\geq 1} c_ng_n(\omega)e^{int}$ with $g_n$ an independent identically
distributed sequence of complex-valued Gaussian random variables, then
$F\in L^p$ if and only if the sequence $(c_n)\in \ell^2$; see
Paley and Zygmund \cite{PZ}.
Also see a recent paper by Ayache and Tzvetkov \cite{aytzve}. A more detailed discussion of the comparison
between classical Besov spaces and Fourier-Besov spaces is given in Subsection \ref{FAR} after the
statement of Theorem \ref{main}. We now turn our attention to the regularity of the Brownian
motion on modulation spaces. By contrast to the classical function spaces discussed above, modulation spaces do capture the
dependence of regularity on one of its defining indexes.

\subsection{Fourier analytic representation}\label{FAR}

We are interested in the local-in-time regularity on modulation spaces $M^{p, q}_s (\T)$ of the periodic
Brownian motion defined on the 1-dimensional torus. The case of the Brownian motion on the $d$-dimensional torus $\T^d$ is similar;
see Remark \ref{Rmk2}. Moreover, we can obtain the regularity of the Brownian motion on the real line as long
as we replace the periodic modulation spaces $M^{p, q}(\T)$ with the localized versions $M^{p, q}_s(I)$,
where $I\subset \R$ is any given bounded interval; see Remark \ref{Rmk1} and Appendix \ref{ApB}.

Without loss of generality, it suffices to consider the mean-zero
complex-valued Brownian loop, that is, satisfying both $\beta(0) =
\beta(2\pi)$ and $\int_0^{2\pi} \beta(t) dt = 0$. The case of a
Brownian loop with non-zero mean follows easily from the
mean-zero case through a translation by the mean of the Brownian
motion. Furthermore, if we let $b(t)$ denote any (non-periodic)
Brownian motion, and set $\beta (t)=b(t)-tb(2\pi)/2\pi$, then $\beta
(t)$ is periodic, and $b(t)$ has the same regularity as $\beta (t)$.

Since, for almost every $\omega\in \Omega$, the Brownian motion
$\beta(t)$ represents a continuous (periodic) function, to simplify the
notation, we will simply write $u$ in the following to denote this
function (for a fixed $\omega$). It is well known \cite{Kahane} that $u$ can be
represented through a {\it Fourier-Wiener series}\footnote{Henceforth,
we  drop a factor of $2 \pi$ when it plays no role.}
as
\begin{equation} \label{2FW}
u(t) =  u(t; \omega)=\sum_{n\not=0}\frac{g_n(\omega)}{n}e^{int},
\end{equation}

\noi where $\{g_n\}$ is a family of independent standard
complex-valued Gaussian random variables, that is,  $\mbox{Re}\, g_n$
and $\mbox{Im}\, g_n$ are independent standard real-valued Gaussian random
variables. Note that we are missing the linear term $g_0(\omega)t$
in the Fourier-Wiener series representation since we are only
considering the (mean zero) periodic case. For the convenience of the
reader, we present a detailed derivation of \eqref{2FW} in Appendix \ref{ApA}.

With the notation above, our result can be stated as follows.
\begin{theorem}\label{main}
Let $1\leq p, q\leq\infty$.
\begin{itemize}
\item[\textup{(a)}] If $q<\infty$, then $u\in M^{p, q}_s (\T)$ a.s. for
$(s-1)q<-1$, and $u\not\in M^{p,q}_s (\T)$ a.s. for $(s-1)q\geq -1$.

\item[\textup{(b)}] If $q=\infty$, then $u\in M^{p, \infty}_s (\T)$ a.s. for
$s<1$, and $u\not\in M^{p, \infty}_s (\T)$ a.s. for $s\geq 1$.

\hspace*{-1.7cm}Moreover,

\item[\textup{(c)}] If $q<\infty$, then $u\in\ft{b}^s_{p,q} (\T)$ a.s. for
$(s-1)p<-1$, and $u\not\in\ft{b}^s_{p,q} (\T)$ a.s. for $(s-1)p\geq -1$.

\item[\textup{(d)}] If $q=\infty$ and $p<\infty$, then
$u\in\ft{b}^s_{p,\infty} (\T)$ a.s. for $(s-1)p\leq-1$, and
$u\not\in\ft{b}^s_{p,\infty} (\T)$ a.s. for $(s-1)p> -1$.
\end{itemize}

\end{theorem}

Note that, in our statements regarding the Fourier-Besov spaces
$\ft{b}^s_{p,q} (\T)$, the case $p=q=\infty$ is already addressed by part
(b), since $\ft{b}^s_{\infty, \infty} (\T)=\mathcal{F} L^{s, \infty} (\T)$. We
think of the case in which $(s-1)p=-1$ as an ``end-point'' that
makes the transition of regularity from $\mathcal{F}L^{s, p} (\T)$ to
$\ft{b}^s_{p, \infty} (\T)$. The case when $(s-1)p=-1$ and $p=2$ (that
is, $s=1/2$) corresponds to the end-point case for the usual Besov
spaces.

It is also worthwhile to note that the Fourier-Besov spaces
$\ft{b}^s_{p, q} (\T)$ are more finely tuned for this analysis than the
regular Besov spaces because of their sensitivity to the value of
$p$ under randomization. Roughly speaking, this is implied by the
definitions of the two spaces: the $\ft{b}^s_{p, q} (\T)$ spaces use the
$L^p_n$ norm on the Fourier side for each dyadic block, while the
$B^s_{p, q}(\T)$ spaces use the $L^p_t$ norm on the physical side for
each dyadic block. Now, in the case of regular Besov spaces, the
$L^p_t$ norms are equivalent under the expectation, that is, the
$L^p_t$ spaces are equivalent under randomization due to a
Khintchine type argument or Paley-Zygmund's theorem; see
\cite{Kahane}. In other words, the $L^p_t$ part is ``insensitive''
to a finite $p$ under randomization.

\smallskip

\begin{remark}\rm\label{Rmk0}
It is well-known that the Gaussian part of the Gibbs measure for the Benjamin-Ono equation corresponds to a
(periodic) fractional Brownian motion having the Fourier-Wiener series representation $\sum_{n\not=0}\frac{g_n (\omega)}{|n|^{1/2}}e^{int}$.
See \cite{TZ3}.
In general, the regularity of any (mean zero periodic) fractional Brownian motion represented by $\sum_{n\not=0}\frac{g_n (\omega)}{|n|^{\alpha}}e^{int}$ is
described by an analogous statement to the one of Theorem \ref{main}. 
(Note that,  when $\al =  0$,  this correspond to the mean zero Gaussian white noise on $\T$.)
In this case, the only modifications are that for the $M^{p, q}_s (\T)$ regularity
we require the condition $(s-\alpha)q<-1$, while for $\ft{b}^s_{p, q} (\T)$ we need $(s-\alpha)p<-1$. These modifications are
straightforward consequences of the summability conditions appearing in the proof of Theorem \ref{main},
and thus  we omit the details.
\end{remark}

\begin{remark}\rm\label{Rmk1}
We can define the local-in-time versions of the functions spaces above in the following way.
Given an interval $I\subset \R$, we let $M^{p, q}_s (I)$ denote the restriction of $M^{p, q}_s(\R)$ onto $I$
via
$$
 \|u\|_{M^{p, q}_s (I)} = \inf \big\{ \|v\|_{M^{p, q}_s (\R)}: v = u \text{ on } I\big\}.
$$

\noi The local-in-time versions of other function spaces can be defined in an analogous manner. With this notation, we can show that, given a
bounded interval $I\subset \R$, \emph{Theorem \ref{main} holds with the same conditions on the indices for the Brownian motion $b(t)$ (defined on
the unbounded domain $\R$) on the spaces $M^{p, q}_s (I)$, $W^{p, q}_s (I)$, and $\ft{b}^s_{p, q} (I)$.} Clearly, $b(t)$ is unbounded in any of
$M^{p, q}_s (\R)$, $W^{p, q}_s (\R)$, or $\ft{b}^s_{p, q} (\R)$ a.s., since $b(t)$ represents a.s.~an unbounded function on $\R$.
Namely, this unboundedness of the norms on $\R$ comes from the integrability condition rather than the differentiability condition. 
The ``equivalence''
of the periodic function spaces (defined on $\T$) and their local-in-time versions (defined on some bounded interval $I$) is proved in detail in
Appendix \ref{ApB}.
\end{remark}

\begin{remark} \rm\label{Rmk2}
\emph{Theorem \ref{main} can also be extended to the Brownian motion on the $d$-dimensional torus $\T^d$.} We simply need to modify the corresponding statements
by imposing the inequality $(s-1)q<-d$ when considering the regularity on $M^{p, q}_s (\T^d)$, or $(s-1)p<-d$ when considering the $\ft{b}^s_{p,q} (\T^d)$
regularity. Let us briefly indicate how these conditions arise. 
%
%Starting with the white noise on
%$\T^d$, $W(t) =  \sum_{n \in \Z^d \setminus \{0\}}g_n (\omega)  e^{in\cdot t},$ and formally taking its antiderivative, we can
%consider 
The periodic Brownian motion on $\T^d$ can be represented by 
\begin{equation}\label{FWZ}
u(t) =
\sum_{n \in \Z^d \setminus \{0\}}\frac{g_n (\omega)}{|n|} \,  e^{in\cdot t},
\end{equation}

\noi
where $\{g_n\}_{n \in \Z^d \setminus \{0\}}$
is a family of independent standard complex-valued Gaussian random variables.
See Bourgain \cite{BO7, BO8}.
Note that \eqref{FWZ} corresponds to a typical element in the support of 
the $d$-dimensional analogue of the Wiener measure $\mu$ defined 
in \eqref{Wiener1} and \eqref{Wiener3}. 
Thus, in the proof of Theorem
\ref{main}, say in \eqref{DIV0}, we only need to change our summation $\sum \jb{n}^{(s-1)q}$ over the lattice $\Z$ (which requires $(s-1)q<-1$ for
convergence) to that over $\Z^d$ to recover the correct condition $(s-1)q<-d$. 
See also Remark \ref{Rmk6}.
We note also that the regularity of the white noise $W(t)$ on
$\T^d$ in \cite{ver1} can be deduced from the $d$-dimensional version of our Theorem \ref{main} by simply dropping one regularity from that
of the Brownian motion $u(t)$. For example, we have that if $q<\infty$, then $W\in\ft{b}^s_{p,q} (\T^d)$ a.s. for
$sp<-d$, and $W\not\in\ft{b}^s_{p,q} (\T^d)$ a.s. for $sp\geq -d$. 
See also \cite{OH4} and Remark \ref{Rmk0}.
Finally, it is worth recalling that Brownian motion in higher
dimensions is commonly referred to as a Gaussian Free Field (GFF). The physics literature uses the names massless free field or Euclidean
bosonic massless free field. For a nice overview of GFF and its applications, see Sheffield's article \cite{shef}.
See also Remark \ref{RmkA}.
\end{remark}

Theorem \ref{main} states that the Brownian motion belongs a.s.~to
$M^{p, q}_s (\T)$, $W^{p, q}_s (\T)$, and $\mathcal{F} L^{s, q}(\T)$ for $(s-1)q <
-1$, and to $\ft{b}^s_{p, q} (\T)$ for $(s-1)p < -1$, and to
$\ft{b}^s_{p, \infty} (\T)$ for $(s-1)p \leq -1$. However, in
applications, it is often very important to know how large the
estimate on the norm is likely to be; see the works by Bourgain
\cite{BO4, BO5, BO7, BO8}, Burq and Tzvetkov \cite{BT1, BT3}, and
the second author \cite{nors, OH3, OH4, OHSBO}. The following theorem
provides us with the desirable ``large deviation estimates''.

\begin{theorem}\label{LDE1}
There exists  $c>0$ such that for (sufficiently large) $K>0$, the
following holds:
\begin{itemize}
\item[\textup{(i)}] If $(s-1)q<-1$, then $Pr\big( \|u(\omega)\|_{M^{p,
q}_s (\T)} > K\big) < e^{-cK^2}.$

\item[\textup{(ii)}] If $(s-1)p<-1$, then $Pr\big(
\|u(\omega)\|_{\ft{b}^s_{p, q} (\T)}> K\big) < e^{-cK^2}.$

\item[\textup{(iii)}] If $(s-1)p=-1$ (and $q = \infty$), then $Pr\big(
\|u(\omega)\|_{\ft{b}^s_{p, \infty} (\T)}> K\big) < e^{-cK^2}.$
\end{itemize}
\end{theorem}

\begin{remark}\rm\label{Rmk4} We note that the same estimates
in Theorem \ref{LDE1} hold for all $K>0$ as long as we replace the
right-hand side by $Ce^{-cK^2}$ for an appropriate $C$. If we wish
to ``normalize'' the constant $C$ so that it equals $1$, we must
select $K>0$ sufficiently large so that $Ce^{-\frac{c}{2}K^2}<1$.
\end{remark}

\noi The proofs of parts (i) and (ii) in Theorem \ref{LDE1} rely on the
theory of abstract Wiener spaces and Fernique's theorem. A detailed
discussion of these proofs and the afferent technicalities is given
in Section 3.

\begin{proof} [Proof of Theorem \ref{main}]

We begin by showing statements (a) and (b). Recall that we have $M^{p, q}_s (\T)=W^{p, q}_s (\T)=\mathcal{F}L^{s, q} (\T)$, where we
defined the latter space via the norm:
\[ \|u\|_{\mathcal{F}L^{s, q}  (\T)} = \| \jb{n}^s \ft{u}(n)\|_{\ell^q_n}.\]

\noi We easily see that, with $\mathbb{E}$ denoting expectation, we have
\begin{equation}\label{DIV0}
\mathbb{E} \big[\| u \|_{\mathcal{F} L^{s, q} (\T)}^q \big] = \sum_{n\ne 0}
\jb{n}^{sq} |n|^{-q}\mathbb{E} \big[ |g_n|^q \big] \sim \sum_{n\ne 0}
\jb{n}^{(s-1)q} < \infty
\end{equation}

\noi if and only if $(s-1) q<-1$.

Let us now define
\begin{equation} \label{XJ}
X^{(q)}_j (\omega) :=2^{-j} \sum_{|n| \sim 2^j} |g_n(\omega) |^q.
\end{equation}

\noi Then, $X_j^{(q)} \to c_q: =\mathbb{E} |g_1|^q $ almost surely,
since $Y^{(q)}_j :=2^{-j} \sum_{1\leq |n| \leq 2^{j-1}} |g_n|^q \to
c_q$ by the strong law of large numbers, and $X^{(q)}_j = 2
Y^{(q)}_{j+1} - Y^{(q)}_j$.

Suppose now that $(s-1) q \geq -1$. Then, we have
\begin{align} \label{DIV1}
\| u \|_{\mathcal{F} L^{s, q} (\T)}^q & = \sum_{n\ne 0} \jb{n}^{sq}|n|^{-q} |g_n
(\omega)|^q
\sim \sum_{j = 0}^\infty \sum_{|n|\sim 2^j}\jb{n}^{(s-1)q} |g_n (\omega)|^q \\
 & \geq \sum_{j = 0}^\infty \sum_{|n|\sim 2^j}\jb{n}^{-1} |g_n (\omega)|^q \notag
\sim \sum_{j = 0}^\infty X^{(q)}_j (\omega) = \infty,  \ \mbox{a.s.}
\end{align}

\noi Hence, when $q < \infty$, $u \in \mathcal{F}L^{s, q} (\T)$ a.s. for
$(s-1)q<-1$, and $u \notin \mathcal{F}L^{s, q} (\T)$ a.s. for $(s-1)q
\geq -1$.

Now, let $q = \infty$. Then, we have
\begin{align*}
\| u \|_{\mathcal{F} L^{s, \infty} (\T)} = \sup_{n \ne 0} \jb{n}^{s}|n|^{-1} |g_n(\omega)|
\sim \sup_{n \ne 0} \jb{n}^{s-1} |g_n(\omega)|< \infty,
\end{align*}

\noi a.s. for $s < 1$ (i.e. ``$(s-1) \cdot \infty < -1$''), since $
\lim_{n \to \infty} n^{-\eps} |g_n (\omega)| = 0$ for any $\eps > 0$. When $s \geq 1$,
the continuity from below of the probability measure gives
\[ Pr( \| u \|_{\mathcal{F} L^{s, \infty} (\T)} < \infty)
\leq Pr(\, \sup_n|g_n (\omega)|< \infty) = \lim_{K\to \infty} Pr(\,
\sup_n|g_n (\omega)|< K) = 0.
\]

\noi Hence, $u\in \mathcal{F}L^{s, \infty} (\T)$ a.s. for $s<1$, and $u
\notin \mathcal{F}L^{s, \infty} (\T)$ a.s. for $s\geq1$.

\medskip

Regarding (c) and (d), recall that we defined $\ft{b}^s_{p, q} (\T)$ via
the norm in \eqref{Besov2} and \eqref{Besov3}.
First, suppose $(s-1) p < -1$. Then, for $1\leq p < \infty$ and
$q \geq 1$, we have
\begin{align} \label{CON1}
\mathbb{E}\big[\|u\|_{\ft{b}^s_{p, q} (\T)}^p\big] &\leq
\mathbb{E}\big[\|u\|_{\ft{b}^s_{p, 1} (\T)}^p\big] \leq \bigg(\sum_{j =
0}^\infty \Big( \sum_{|n|\sim 2^j}  \jb{n}^{sp}|n|^{-p}
\mathbb{E}[|g_n|^p] \Big)^\frac{1}{p}\bigg)^p\\
&\sim \bigg(\sum_{j = 0}^\infty \Big(\sum_{|n|\sim 2^j}
\jb{n}^{(s-1)p}\Big)^\frac{1}{p}\bigg)^p \sim \Big(\sum_{j =
0}^\infty 2^{\frac{(s-1)p+1}{p}j}\Big)^p < \infty. \notag
\end{align}

\noi Also, when $p, q < \infty$, $\|u\|_{\ft{b}^s_{p, q} (\T)} = \infty$
a.s. for $(s-1)p \geq -1$ by modifying the argument in \eqref{DIV1}.
When $p = \infty$ and $s  < 1$, we have $\|u\|_{\ft{b}^s_{\infty,
q} (\T)} \leq \|u\|_{\ft{b}^s_{r, q} (\T)}$ for $r < \infty$.
Moreover, we can take $r$ to be sufficiently large
such that $(s-1) r <-1$. Then, it follows from \eqref{CON1} that
$\|u\|_{\ft{b}^s_{\infty, q} (\T)} < \infty$ a.s. For $s \geq 1$, we have
\begin{align*}
\|u\|_{\ft{b}^s_{\infty, q} (\T)} \sim  \big\| \sup_{|n|\sim 2^j} \jb{n}^{s-1}|g_n|
\big\|_{\ell^q_j} \geq \Big\| 2^{-j} \sum_{|n|\sim 2^j} |g_n|
\Big\|_{\ell^q_j} = \infty,  \  \mbox{a.s.},
\end{align*}

\noi as in \eqref{DIV1}. Hence,  $u \in \ft{b}^s_{p, q} (\T)$ a.s. for
$(s-1)p<-1$, while $u \notin \ft{b}^s_{p, q} (\T)$ a.s.
for $(s-1)p \geq -1$ and $q < \infty$.

Finally, we consider the case $ q = \infty$ and $p < \infty$. In the
endpoint case $(s-1) p = -1$, we have
\begin{align} \label{CON2}
\|u\|_{\ft{b}^s_{p, \infty} (\T)}^p
& \sim \sup_j  \sum_{|n| \sim 2^j} \jb{n}^{(s-1)p} |g_n (\omega)|^p \\
& \sim \sup_j \, 2^{-j} \sum_{|n| \sim 2^j} |g_n (\omega)|^p = \sup_{j}
X^{(p)}_j (\omega)< \infty, \ \text{a.s.},  \notag
\end{align}

\noi where $X^{(p)}_j$ is defined in \eqref{XJ}. When $(s-1) p >
-1$, a similar computation along with the convergence of $X^{(p)}_j$
shows that $u \notin \ft{b}^s_{p, \infty} (\T)$ a.s.
\end{proof}

\subsection{Alternate proof for the Besov spaces}

We close this section by providing a new (alternate) proof, via
Fourier-Wiener series \eqref{2FW}, of the regularity results on
Besov spaces that we exposed at the end of Subsection 2.1. We
decided to include this proof because random Fourier series are a
unifying theme of this paper. Furthermore, our proof, which is done
on the Fourier side, seems to complement nicely the existing one
using Franklin-Wiener series on the physical side \cite{ROY}. For notational
simplicity we write $B^s_{p, q}$ for $B^s_{p, q} (\T)$.

We begin by recalling the general Gaussian bound
\begin{equation}\label{Gaussbd}
 \Big\| \sum_n c_n g_n(\omega) \Big\|_{L^p(\Omega)} \leq C \sqrt{p}
 \|c_n\|_{l^2_n};
\end{equation}

\noi see the works of Burq and Tzvetkov \cite{BT2}, and Tzvetkov
\cite{TZ3}. Also, see Lemma \ref{PROP:hyp} below.
Then, we have, for $1 \leq p < \infty$ and $1\leq q \leq
\infty$, (recall $t \in \T = \R /2\pi \Z$),
\begin{align*}
\mathbb{E} \big[\|u\|_{B^s_{p, q}}\big]
& = \mathbb{E}\bigg\| \Big\| \sum_{|n|\sim 2^j} \jb{n}^{s}|n|^{-1} g_n
e^{int} \Big\|_{L^p_t} \bigg\|_{\ell^q_j}
\leq \bigg\| \Big\| \big\|\sum_{|n|\sim 2^j} \jb{n}^{s}|n|^{-1}
 g_n(\omega) e^{int}\big\|_{L^p(\Omega)} \Big\|_{L^p_t} \bigg\|_{\ell^1_j}\\
& \lesssim \sum_{j = 0}^\infty \Big( \sum_{|n|\sim 2^j
}\jb{n}^{2(s-1)} \Big)^\frac{1}{2} < \infty
\end{align*}

\noi for $2(s-1) < -1$, i.e. $s < \frac{1}{2}$. When $s <
\frac{1}{2}$ and $p = \infty$, Sobolev's inequality gives
$\|\jb{\dt}^s u \|_{L^\infty} \lesssim \|\jb{\dt}^{s + \eps} u
\|_{L^r}$, for small $\eps >0$ and large $r$ such that $ s + \eps <
\frac{1}{2}$ and $\eps r >1$. Then, the above computation shows that
$\|u\|_{B^s_{\infty, q}} \lesssim \|u\|_{B^{s+\eps}_{r, q}} <
\infty$ a.s. for $ s < \frac{1}{2}$.

Suppose $2 \leq p\leq \infty$. Then, for $q < \infty$ and $s \geq
\frac{1}{2}$, we have
\begin{align} \label{DIV2}
\|u\|_{B^s_{p, q}} \geq \|u\|_{B^s_{2, q}} \sim \Big\|
\Big(\sum_{|n|\sim 2^j} 2^{2(s-1)j} |g_n|^2   \Big)^\frac{1}{2}
\Big\|_{\ell^q_j} \sim \big\| \big(2^{(2s-1)j} X^{(2)}_j
\big)^\frac{1}{2} \big\|_{\ell^q_j} = \infty,  \ \text{a.s.}
\end{align}

\noi since $X^{(2)}_j$ defined in \eqref{XJ} converges to $c_2 > 0$
a.s. It also follows from \eqref{DIV2} that $\|u\|_{B^s_{p, \infty}}
= \infty$ a.s. for $s > \frac{1}{2}$ when $ q = \infty$. Now, let
$1\leq p < 2$, $q < \infty$, and $s \geq \frac{1}{2}$. Then, we have
\begin{align*}
\|u\|_{B^s_{p, q}}^q \geq \|u\|_{B^\frac{1}{2}_{1, q}}^q \sim \sum_j
Z_j^{(q)} (\omega)  ,
\end{align*}

\noi where $Z_j^{(q)}(\omega) = \|  \wt{X}_j (t, \omega) \|_{L^1_t
}^q$ and
\begin{align} \label{XJ2}
\wt{X}_j (t, \omega) & :=2^{-\frac{j}{2}} \sum_{|n| \sim 2^j}
g_n(\omega)e^{int}.
\end{align}

\noi Note that for each $t \in \T$, $\wt{X}_j (t, \omega)$ is a
standard complex Gaussian random variable. Thus, we have
$\mathbb{E}\big[ \|\wt{X}_j (t, \omega)\|_{L^1_t}\big] =
\big\|\mathbb{E} \,|\wt{X}_j (t, \omega)|\big\|_{L^1_t} = 2\pi c_1
>0$. In particular, we have
\begin{equation} \label{Gaussbd2}
\big( \mathbb{E} \big[Z_j^{(q)}\big]\big)^\frac{1}{q}
\geq \mathbb{E}\big[ \|\wt{X}_j (t, \omega)\|_{L^1_t}\big] = 2\pi c_1
\end{equation}

\noi Also, by \eqref{Gaussbd}, we have $\mathbb{E}\big[
|Z_j^{(q)}|^2\big] \leq \big\| \|\wt{X}_j (t,
\omega)\|_{L^{2q}_\omega} \big\|^{2q}_{L^1_t} \leq C (2q)^q <\infty. $
Hence, by Kolmogorov's strong law of large numbers, we have
$\frac{S_m- \mathbb{E}[S_m]}{m+1} \to 0$ a.s., where $ S_m =
\sum_{j = 0}^m Z_j^{(q)}$. It follows from \eqref{Gaussbd2} that
$\frac{\mathbb{E}[S_m]}{m+1} \geq (2\pi c_1)^q.$ This
implies that $S_m (\omega)\to \infty $ a.s.~in $\omega\in\Omega$.
Hence, we have
\[ \|u\|_{B^s_{p, q}}
\geq \|u\|_{B^\frac{1}{2}_{1, q}} = \Big(\lim_{m\to \infty} S_m
\Big)^\frac{1}{q} = \infty, \ \text{a.s.}\]

For $q = \infty$ and $ s > \frac{1}{2}$, we have
\begin{align} \label{Besovlowerbd1}
\|u\|_{B^s_{p, \infty}} \geq \|u\|_{B^s_{1, \infty}} \sim \sup_j
2^{(s-\frac{1}{2})j}Z_j^{(1)} (\omega)  ,
\end{align}

\noi where $Z_j^{(1)}(\omega) = \|  \wt{X}_j (t, \omega) \|_{L^1_t
}$. Note that $\mathbb{E}\big[Z_j^{(1)}\big] = 2 \pi c_1$ and
$\mathbb{E}\big[|Z_j^{(1)}|^2\big] \leq C < \infty$ for all $j$.
This implies that there exist $\dl, \eps > 0$ and $\Omega_j$, $j =
0, 1, \dots$, such that $Z_j^{(1)}(\omega) > \dl $ for $\omega \in
\Omega_j$ and $Pr (\Omega_j) > \eps$. In particular, we have
$\sum_j\, Pr (\Omega_j) = \infty$. Then, by the Borel zero-one law,
\begin{equation} \label{Besovlowerbd2}
 Pr ([Z_j^{(1)}(\omega) > \dl, \text{ infinitely often}]) = 1.
\end{equation}

\noi From \eqref{Besovlowerbd1} and \eqref{Besovlowerbd2},
$\|u\|_{B^s_{p, \infty}} = \infty$ a.s. for $s > \frac{1}{2}$.

Now, consider $p = q = \infty$ and $ s = \frac{1}{2}$. We have
\[ \|u\|_{B^\frac{1}{2}_{\infty, \infty}} = \sup_j \| \wt{X}_j(t, \omega) \|_{L^\infty_t}
\geq \sup_j | \wt{X}_j(t_j^*, \omega) |,
\]

\noi where $\wt{X}_j$ is defined in \eqref{XJ2} and $t^*_j$'s are
points in $\T$. Recall that $\{\wt{X}_j(t_j^*)\}_{j = 0}^\infty$ is
a family of independent standard complex-valued Gaussian random
variables. Hence, $\sup_j | \wt{X}_j(t_j^*, \omega) | = \infty$ a.s.
and thus $\|u\|_{B^\frac{1}{2}_{\infty, \infty}} = \infty$ a.s.

Finally, we consider the case $p < \infty$, $q = \infty$ and $s =
\frac{1}{2}.$ First, assume $p \leq 2$. Then, we have
\begin{align*}
\| u\|^2_{B^\frac{1}{2}_{p, \infty}} \leq \|
u\|^2_{B^\frac{1}{2}_{2, \infty}} = \sup_j X^{(2)}_j < \infty, \
\text{a.s.}
\end{align*}

\noi since $X^{(2)}_j$ defined in \eqref{XJ} converges to $c_2$ a.s.
In the following, we consider
\begin{align} \label{Bes0}
\| u\|_{B^\frac{1}{2}_{p, \infty}}^p \sim\sup_j \Big\|2^{-\frac{j}{2}}
\sum_{|n|\sim 2^j} g_n e^{int}\Big\|_{L^p_t}^p
\end{align}

\noi only for $p = 2k$ with $k = 2, 3, \dots$ since $\|
u\|_{B^\frac{1}{2}_{p, \infty}} \leq \| u\|_{B^\frac{1}{2}_{2k,
\infty}}$ for $p \leq 2k$.

When $p = 4$, we have
\begin{align} \notag
\Big\|2^{-\frac{j}{2}} \sum_{|n|\sim 2^j} g_n
e^{int}\Big\|_{L^4_t}^4 = \ & 2\cdot 2^{-2j} \sum_{|n_1|, |n_2| \sim 2^j}
|g_{n_1}|^2|g_{n_2}|^2
 + 2^{-2j} \sum_{\substack{|n_\al|, |m_{\beta}| \sim 2^j \\
 n_\al \ne m_{\beta}\\n_1 + n_2 = m_1 + m_2}} g_{n_1}g_{n_2}\cj{g}_{m_1}\cj{g}_{m_2}\\
& - 2^{-2j} \sum_{|n| \sim 2^j} |g_{n}|^4 =: \I^{(2)}_j +
\II^{(2)}_j + \III^{(2)}_j. \label{B4}
\end{align}

\noi Note that $\I^{(2)}_j = \big(X_j^{(2)} \big)^2\to c_2^2$ and
$\III^{(2)}_j = 2^{-j} X_j^{(4)}\to 0$ a.s. by the strong law of
large numbers, where $X_j^{(p)}$ is defined in \eqref{XJ}. Hence, it
suffices to prove that $ \sup_j |\II^{(2)}_j| < \infty$ a.s. By
independence of $\II^{(2)}_j$ and the Borel zero-one law, it suffices to
show that
\begin{equation} \label{Bes1}
\sum_j Pr( |\II^{(2)}_j| > K) < \infty
\end{equation}

\noi for some $K>0$. By Chebyshev's inequality, we have
\begin{equation} \label{cheby1}
 Pr( |\II^{(2)}_j| > K) \leq K^{-2} \mathbb{E} \big[|\II^{(2)}_j|^2 \big]
\leq C_2 K^{-2} 2^{-j}.
\end{equation}

\noi Hence, \eqref{Bes1} follows, and thus we have $\|
u\|_{B^\frac{1}{2}_{4, \infty}} <\infty$ a.s.

In order to estimate \eqref{Bes0} for the general case $p = 2k$, we
use an induction argument and assume the existence of estimates for
$p = 2, \dots, 2(k-1)$. We have
\begin{align} \notag
\Big\|2^{-\frac{j}{2}} \sum_{|n|\sim 2^j} g_n
e^{int}\Big\|_{L^{2k}_t}^{2k} = \ & k!\cdot 2^{-kj} \sum_{|n_\al| \sim 2^j}
\prod_{j = 1}^k|g_{n_\al}|^2
 + 2^{-kj} \sum_{*}\prod_{\al = 1}^k g_{n_\al}\prod_{\beta = 1}^k\cj{g}_{m_{\beta}}\\
& + \text{ error terms} =: \I^{(k)}_j + \II^{(k)}_j + \text{ error
terms}, \label{BP}
\end{align}

\noi where $* = \{n_\al, m_\beta : \al, \beta = 1, \dots, k, \,
|n_\al|, |m_{\beta}| \sim 2^j, \ n_\al \ne m_{\beta}\ \sum n_\al =
\sum m_{\beta}\}$. Note that $\I^{(k)}_j$ consists of the terms for
which $n_\al$'s and $m_{\beta}$'s form exactly $k$ pairs (including
higher multiplicity) and that $\II^{(k)}_j$ consists of the terms
with no pair. There are two types of error terms, which we call of
type (i) and type (ii):
\begin{itemize}
\item[(i)]$\text{error}^{(k)}_j\text{(i)}$:  $\{ n_\al\}_{\al = 1}^k = \{ m_{\beta}\}_{\beta = 1}^k$,
and there exists at least one pair $\al$, $\wt{\al}$ with $\al \ne
\wt{\al}$ such that $n_\al = n_{\wt{\al}}$, i.e., $n_\al$'s and
$m_{\beta}$'s form exactly $k$ pairs, but there exists at least ``four
of a kind'' (or higher order): $n_{\al_1} = n_{\al_2}= m_{\beta_1}=
m_{\beta_2}$ with $\al_1 \ne \al_2$ and $\beta_1 \ne \beta_2$.

\item[(ii)]$\text{error}^{(k)}_j\text{(ii)}$: $\{ n_\al\}_{\al = 1}^k \ne \{ m_{\beta}\}_{\beta = 1}^k$,
but there exists at least one pair $\al$, $\beta$ such that $n_\al =
m_{\beta}$. i.e. $n_\al$'s and $m_{\beta}$'s form exactly $\ell$ pairs
for some $ 1\leq \ell \leq k - 1$.
\end{itemize}

\noi As before, by the strong law of large numbers, $\I^{(k)}_j =
\big(X_j^{(2)} \big)^k\to c_2^k$ a.s.~and the error terms of type
(i) go to 0 a.s.~since each of them can be written as
\begin{align*}
2^{-(k-L)j} \prod_{\ell \, : \,  k = \sum_{\ell = 1}^L k_\ell} X_j^{(2k_\ell)},
\ \ L < k,
\end{align*}

\noi where each $X_j^{(2k_\ell)} \to c_{2k_\ell}$ a.s. Hence, we have
$Pr(\sup_j |\text{error}^{(k)}_j\text{(i)}| <\infty) = 1$. As for
the error terms of type (ii), the worst ones can be written as
\begin{align*}
2^{-kj} \sum_{n} |g_n|^2 \sum_{**}\prod_{\al = 1}^{k-1}
g_{n_\al}\prod_{\beta = 1}^{k-1}\cj{g}_{m_{\beta}}
\end{align*}

\noi where $** = \{n_\al, m_\beta : \al, \beta  = 1, \dots, k-1, \,
|n_\al|, |m_{\beta}| \sim 2^j, \ n_\al \ne m_{\beta}\ \sum n_\al =
\sum m_{\beta}\}$. It is basically a product of $X^{(2)}_j$ (which
converges to $c_2$ a.s.) and $\II^{(k-1)}_j$, i.e., the $(k-1)$-fold
products over frequencies $\{n_\al\}_{\al = 1}^{k-1}$ and $ \{
m_{\beta}\}_{\beta = 1}^{k-1}$ containing no pair, which appeared at
the $k-1$ inductive step. All the other error terms of type (ii) can
be basically written as $\text{error}^{(\ell)}_j\text{(i)} \cdot
\II^{(k-\ell)}_j$ for some $\ell = 1, \dots, k-1$. Hence, we have
$Pr(\sup_j |\text{error}^{(k)}_j\text{(ii)}| <\infty) = 1$.

Now, it remains to estimate $\II^{(k)}_j$. As before, it suffices to
show that
\begin{equation} \label{Bes2}
\sum_j Pr( |\II^{(k)}_j| > K) < \infty
\end{equation}

\noi for some $K>0$. By Chebyshev's inequality, we have
\[ Pr( |\II^{(k)}_j| > K) \leq K^{-2} \mathbb{E} \big[|\II^{(k)}_j|^2 \big]
\leq C_k K^{-2} 2^{-j}.\]

\noi Hence, \eqref{Bes2} follows, and thus we have $\|
u\|_{B^\frac{1}{2}_{2k, \infty}} <\infty$ a.s. This completes the
proof.

\section{Large deviation estimates}\label{SLDE}

\subsection{Abstract Wiener spaces and Fernique's theorem}\label{AWS}

Let $B$ denote any of the spaces $M^{p, q}_s (\T)$, $W^{p, q}_s (T)$,
$\mathcal{F} L^{s, q} (\T)$, or $\ft{b}^s_{p, q} (\T)$, and, as before, let
$u$ be the mean zero complex-valued Brownian loop on $\T$. While the
previous section was concerned with the question of $B$-regularity,
i.e., whether or not $u$ is in $B$, this section will be concerned
with the complementary topic of large deviation estimates on $B$.
Specifically, we will establish estimates of the form
\begin{equation} \label{L}
Pr\big( \|u(\omega)\|_B > K\big) < e^{-cK^2},
\end{equation}
for large $K>0$, and some constant $c=c(B)>0$. As we shall see, the
theory of abstract Wiener spaces and Fernique's theorem play a
crucial role in establishing estimates such as \eqref{L} on all
non-endpoint cases, see Proposition \ref{PROP:FER} below. In
Subsection 3.2, we prove that the large deviation estimate still
holds for $\ft{b}^s_{p, \infty} (\T)$ in the endpoint case $(s-1)p = -1$
even though Fernique's theorem is not applicable. We also discuss,
in Subsection 3.3, the issue of deviation estimates in the endpoint
case of regular Besov spaces $B_{p, \infty}^{\frac{1}{2}}, 1\leq p<\infty.$ An
alternate point of view can be found in \cite{hytver}. For non-endpoint deviation estimates on
$B_{p, q}^{s}$, the reader is referred to Roynette's work \cite{ROY}.

Recall that if $u$ is the mean zero complex-valued Brownian loop on
$\T$, then we can expand it in its Fourier-Wiener series as
\begin{equation} \label{FW1}
u(x,\omega) = \sum_{n\ne 0} \frac{g_n(\omega)}{n} e^{inx}, \quad x
\in \T,
\end{equation}
\noi where $\{g_n(\omega)\}_{n \ne 0}$ is a family of independent
standard complex-valued Gaussian random variables. This induces a
probability measure on the periodic functions on $\T$, namely the
mean zero\,{\it Wiener measure} on $\T$, which can be formally written
as
\begin{equation} \label{Wiener1}
d \mu = Z^{-1} \exp\bigg(-\frac{1}{2}\int_\T |u_x|^2 dx \bigg)
\prod_{x \in \T} d u(x), \ \ u \text{ mean } 0.
\end{equation}
 
 \noi
 In the following, we use the theory of abstract Wiener spaces to
provide the precise meaning of expression \eqref{Wiener1}. Let $u(x)
= \sum_{n\ne 0} \ft{u}_n e^{inx}$ denote any periodic function on
$\mathbb{T}$ with mean 0. We define the Gaussian measure $\mu_N$ on
$\mathbb{C}^{2N}$  with the  density
\begin{equation}\label{Wiener2}
 d \mu_N = Z_N^{-1} \exp\bigg(-\frac{1}{2}\sum_{0<|n|\leq N}|n|^2  |\ft{u}_n|^2 \bigg)
 \prod_{ 0<|n|\leq N } d \ft{u}_n ,
\end{equation}

\noindent where $d \ft{u}_n$ denotes the complex Lebesgue measure on
$\C$ and
\[Z_{N} = \int_{\mathbb{C}^{2N}} \exp\bigg(-\frac{1}{2}
\sum_{0<|n|\leq N} |n|^2 |\ft{u}_n|^2\bigg) \prod_{0<|n|\leq N}  d
\ft{u}_n . \]

In our definition above, we have abused the notation and we denoted
{\it any} generic periodic function on $\T$ (not just the Brownian
motion) by the letter $u$. The context, however, will make it clear
when we refer specifically to Brownian motion.

Note that the measure $\mu_N$ is the induced probability measure on
$\mathbb{C}^{2N}$ (that is, the $2N$ dimensional complex Gaussian or
$4N$ dimensional real Gaussian measure) under the map $ \omega
\mapsto \{ g_n(\omega)/|n| \}_{0< |n|\leq N}$. Indeed, if we replace
$\ft{u}_n$ by $g_n/|n|$ in \eqref{Wiener2}, we have
$$d \mu_N = \wt{Z_N}^{-1}\prod_{ 0<|n|\leq N}\exp\bigg(-\frac{1}{2}|g_n|^2\bigg)\, dg_n,$$
where
$$\wt{Z_N}=\prod_{ 0<|n|\leq N}\int_{\mathbb{C}}\exp\bigg(-\frac{1}{2}|g_n|^2\bigg)\, dg_n=(2\pi)^{2N}.$$

We would like to define the mean zero Wiener measure in
\eqref{Wiener1} as a limit of the finite dimensional Gaussian
measures $\mu_N$ as $N\to \infty$, i.e., we would like to define the
Wiener measure $\mu$ in \eqref{Wiener1} by
\begin{equation} \label{Wiener3}
 d \mu = Z^{-1} \exp\bigg(-\frac{1}{2}\sum_{n\ne 0 } |n|^2 |\ft{u}_n|^2\bigg) \prod_{ n\ne 0} d \ft{u}_n ,
\end{equation}
where
$$Z = \int \exp\bigg(-\frac{1}{2}\sum_{n\ne 0 } |n|^2
|\ft{u}_n|^2\bigg) \prod_{  n \ne 0 } d \ft{u}_n.$$ Note that the
expression in the exponent in \eqref{Wiener3} can be written as
\[-\frac{1}{2}
\sum_{n\ne 0 } |n|^2 |\ft{u}_n|^2 = -\frac{1}{2}\|u\|_{\dot{H}^1} =
-\frac{1}{2} \Big\langle |\dx|^{2-2s} |\dx|^s u, |\dx|^s
u\Big\rangle_{L^2} = -\frac{1}{2} \jb{ B_s^{-1} u, u }_{\dot{H}^s},
\]

\noi where $B_s = |\dx|^{2s-2}$.

It follows from the theory of Gaussian measures on Hilbert spaces
that \eqref{Wiener3} defines a countably additive measure on
$\dot{H}^s$ if and only if $B_s$ is of trace class, i.e., if
$\sum_{n \ne 0} |n|^{2s-2} < \infty$, which is  equivalent to $s
<\frac{1}{2}$; see Zhidkov's work \cite{Z}. This makes the Sobolev space
$H^s$, $s <\frac{1}{2}$, a strong and natural candidate for the
study of Brownian motion. Unfortunately, the spaces under
consideration are not Hilbert spaces in general. To deal with this
issue, the concept of abstract Wiener space comes to the rescue,
since, roughly speaking, it provides us with a larger (Hilbert
or Banach) space, as an extension of $\dot{H}^1$, on which $\mu$ can
be realized as a countably additive probability measure.

In the following, we recall first some basic definitions from Kuo's
monograph \cite{KUO}. Given  a real separable Hilbert space $H$ with
norm $\|\cdot \|_H$, let $\mathcal{F}$ denote the set of finite
dimensional orthogonal projections $\mathbb{P}$ of $H$. Then, a {\it
cylinder set} $E$ is defined by  $E = \{ u \in H: \mathbb{P}u \in
F\}$ where $\mathbb{P} \in \mathcal{F}$ and $F$ is a Borel subset of
$\mathbb{P}H$. We let $\mathcal{R} $ denote the collection of all
such cylinder sets. Note that $\mathcal{R}$ is a field but not a
$\s$-field. Then, the Gaussian measure $\mu$ on $H$ is defined by
\[
 \mu(E) = (2\pi)^{-\frac{n}{2}} \int_F e^{-\frac{1}{2}\|u\|_H^2} du  \]

\noindent for $E \in \mathcal{R}$, where $n = \text{dim} \,
\mathbb{P} H$ and $du$ is the Lebesgue measure on $\mathbb{P}H$. It
is known that $\mu$ is finitely additive but not countably additive
on $\mathcal{R}$.

A seminorm $|||\cdot|||$ in $H$ is called {\it measurable} if, for
every $\eps>0$, there exists $\mathbb{P}_\eps \in \mathcal{F}$ such
that
\begin{equation}\label{msble}
 \mu( ||| \mathbb{P} u ||| > \eps  )< \eps
\end{equation}

\noindent for $\mathbb{P} \in \mathcal{F}$ orthogonal to
$\mathbb{P}_\eps$. Any measurable seminorm  is weaker  than the norm of
$H$, and $H$ is not complete with respect to $|||\cdot|||$ unless
$H$ is finite dimensional. Let $B$ be the completion of $H$ with
respect to $|||\cdot|||$ and denote by $i$ the inclusion map of $H$
into $B$. The triple $(i, H, B)$ is called an {\it abstract Wiener
space}. (The pair $(B,\mu)$ is often called an  abstract Wiener
space as well.)

Now, regarding $v \in B^\ast$ as an element of $H^\ast \equiv H$ by
restriction, we embed $B^\ast $ in $H$. Define the extension of
$\mu$ onto $B$ (which we still denote by $\mu$) as follows. For a
Borel set $F \subset \R^n$, set
\[ \mu ( \{u \in B: ((u, v_1), \cdots, (u, v_n) )\in F\})
:= \mu ( \{u \in H: (\jb{u, v_1}_H, \cdots, \jb{u, v_n}_H )\in
F\}),\]

\noindent where $v_j$'s are in $B^\ast$ and $(\cdot , \cdot )$
denote the dual pairing between $B$ and $B^\ast$. Let
$\mathcal{R}_B$ denotes the collection of cylinder sets $ \{u \in
B:((u, v_1), \cdots, (u, v_n) )\in F \}$ in $B$.

\begin{proposition}[Gross \cite{GROSS}]
$\mu $ is countably additive in the $\s$-field generated by
$\mathcal{R}_B$.
\end{proposition}

In the context of our paper, let $H= \dot{H}^1(\mathbb{T})$. Then,
we have

\begin{theorem} \label{THM:meas}
The seminorms  $\|\cdot\|_{M^{p, q}_s (\T)}$, $\|\cdot\|_{W^{p, q}_s (\T)}$,
and $\|\cdot\|_{\mathcal{F}L^{s, q} (\T)}$ are measurable for $(s-1)q <
-1$. Also, the seminorm $\|\cdot\|_{\ft{b}^s_{p, q} (\T)}$ is measurable
for $(s-1)p < -1$.
\end{theorem}

\noi The proof of Theorem \ref{THM:meas} follows closely the ideas
from \cite[Proposition 3.4]{OH4}. For completeness, we present it in
detail at the end of this subsection.

\begin{corollary}\label{cor:meas}
Let $\mu$ be the mean zero Wiener measure on $\T$. Then, $(M^{p,
q}_s (\T),\mu)$,  $(W^{p, q}_s (\T), \mu)$, and $(\mathcal{F}L^{s, q} (\T), \mu)$
are abstract Wiener spaces for $(s-1)q < -1$. Also, $(\ft{b}^s_{p,
q} (\T),\mu)$ is an abstract Wiener space for $(s-1)p < -1$.
\end{corollary}

\smallskip

\begin{remark} \rm \label{Rmk5} As we shall see later, condition \eqref{msble} is not satisfied for
the endpoint case $\ft{b}^s_{p, \infty} (\T)$ with $(s-1)p = -1$.
Nevertheless, we can still establish a large deviation estimate
using a different approach.
\end{remark}

\begin{remark}\rm \label{Rmk6} 
By making an analogous argument, we can define the mean zero Wiener measure on $\T^d$,
where a typical element on the support is represented by \eqref{FWZ}.
Then, Theorem \ref{THM:meas} and Corollary \ref{cor:meas} can be extended to
$\T^d$. One simply needs to modify the conditions on the indices to read either $(s-1)q<-d$ or
$(s-1)p<-d$ depending on the space considered, that is, modulation or Fourier-Besov, respectively.
\end{remark}

Given an abstract Wiener space $(B, \mu)$, we have the following
integrability result due to Fernique \cite{FER}.

\begin{proposition} [Theorem 3.1 in \cite{KUO}]
\label{PROP:FER} Let $(B,\mu)$ be an abstract Wiener space. Then,
there exists $ c > 0$ such that $ \int_B e^{c \|u\|_B^2} \mu(d u) <
\infty$. In particular, this implies the following large deviation
estimate: there exists $ c' > 0$ such that
\begin{equation} \label{Large1}
\mu(\|u\|_B \geq K ) \leq e^{-c'K^2},
\end{equation}
for sufficiently large $K>0$.
\end{proposition}

From Theorem \ref{THM:meas} and Proposition \ref{PROP:FER}, we
obtain the following corollary.
\begin{corollary}
Let $\mu$ be the mean zero Wiener measure on $\T$. Then, the large
deviation estimate \eqref{Large1} holds for $B = M^{p, q}_s (\T)$, $W^{p,
q}_s (\T)$, and $\mathcal{F}L^{s, q} (\T)$ with $(s-1)q < -1$, Also,
\eqref{Large1} holds for $B = \ft{b}^s_{p, q} (\T)$ with $(s-1)p < -1$.
\end{corollary}

While Proposition \ref{PROP:FER} is not applicable to the endpoint
case $\ft{b}^s_{p, \infty} (\T)$ with $(s-1)p = -1$ (see Remark \ref{Rmk5}), we can
still prove the following result.

\begin{theorem} \label{THM:end2}
Let $(s-1)p = -1$. Then,
\begin{equation} \label{Large2}
\mu(\|u\|_{\ft{b}^s_{p, \infty} (\T)} \geq K ) \leq e^{-cK^2},
\end{equation}
for sufficiently large $K>0$.
\end{theorem}

\noi We prove Theorem \ref{THM:end2} in Subsection
\ref{SUBSEC:Large}. Theorem \ref{THM:end2} also holds for the
endpoint case of the usual Besov spaces $B^s_{p, \infty}$, with $s =
\frac{1}{2}$ and $p<\infty$. However, the proof becomes rather
cumbersome for large values of $p$. Therefore, we will only map the
proof of the large deviation estimates for the Besov spaces with
$p\leq 4$; see Subsection 3.3.

For the proof of Theorem \ref{THM:meas}, we will need the following
lemma from \cite{OHSBO}, which we now recall.

\begin{lemma}[Lemma 4.7 in \cite{OHSBO}] \label{CL:decay}
Let $\{g_n\}$ be a sequence of independent standard complex-valued
Gaussian random variables. Then, for $M$ dyadic and $\dl <
\frac{1}{2}$, we have
\[ \lim_{M\to \infty} M^{2\dl} \frac{\max_{|n|\sim M } |g_n|^2}{ \sum_{|n|\sim M} |g_n|^2} = 0\, \text{ a.s.} \]
\end{lemma}

With these preliminaries, we are ready to prove Theorem
\ref{THM:meas}.

\begin{proof}[Proof of Theorem \ref{THM:meas}]
First, note that $\|u\|_{M^{p, q}_s (\T)}= \|u\|_{W^{p, q}_s (\T)}
=\|u\|_{\mathcal{F}L^{s, q} (\T)} = \|u\|_{\ft{b}^s_{q, q} (\T)}$. Hence, it
suffices to prove the result for $\ft{b}^s_{p, q} (\T)$ with $(s-1)p <
-1$ and any $q \in [1, \infty]$. In view of \eqref{msble},  it
suffices to show that for given $\eps> 0$, there exists large $M_0$
such that
\begin{equation} \label{msble2}
 \mu \big(\|\mathbb{P}_{>M_0}u\|_{\ft{b}^s_{p, q} (\T)} > \eps) < \eps,
\end{equation}

\noi where $\mathbb{P}_{>M_0}$ is the Dirichlet projection onto the frequencies
$|n| > M_0$.

Since $\ft{b}^s_{p, 1} (\T)\subset \ft{b}^s_{p, q} (\T)$, it suffices to prove
\eqref{msble2} for $q = 1$. If $p < 2$ with $(s-1) p < -1$, then by
H\"older inequality, we have
\begin{align*}
\|\jb{n}^s \ft{u}(n)\|_{\ell^p_{|n|\sim2^j}} & \leq
\|\jb{n}^{-\frac{2-p}{2p}}\|_{\ell^{\frac{2p}{2-p}}_{|n|\sim2^j}}
\|\jb{n}^{s + \frac{2-p}{2p}} \ft{u}(n)\|_{\ell^2_{|n|\sim2^j}}\\
& \sim \|\jb{n}^{s + \frac{2-p}{2p}} \ft{u}(n)\|_{\ell^2_{|n|\sim2^j}},
\end{align*}

\noi where $(s + \frac{2-p}{2p} -1)\cdot 2 < -1$. Hence, it suffices
to prove \eqref{msble2} for  $2\leq p \leq \infty$ and $ q = 1$ with
$(s-1) p < -1$.

Now, we consider the case $2 \leq p < \infty$. In the following, we
assume that $u$ is of the form \eqref{FW1}. Fix $K> 1$ and $\dl \in
( 0, \frac{1}{2})$ (to be chosen later.) Then, by Lemma
\ref{CL:decay} and Egoroff's theorem, there exists a set $E$ such
that $\mu (E^c) < \frac{1}{2}\eps$ and the convergence in Lemma
\ref{CL:decay} is uniform on $E$, i.e. we can choose dyadic $M_0$
large enough such that
\begin{equation} \label{A:decay}
\frac{\| \{g_n (\omega) \}_{|n| \sim M} \|_{\ell^{\infty}_n} } {\|
\{g_n (\omega) \}_{|n| \sim M} \|_{\ell^{2}_n} }
 \leq M^{-\dl},
\end{equation}

\noindent for  all $\omega \in E$  and dyadic $M > M_0$. In the
following, we will work only on $E$ and  drop `$\cap E$' for
notational simplicity. However, it should be understood that all the
events are under the intersection with $E$ so that  \eqref{A:decay}
holds.

Let $\{\s_j \}_{j \geq 1}$ be a sequence of positive numbers such
that $\sum \s_j = 1$, and let $ M_j = M_0 2^j$ dyadic. Note that
$\s_j = C 2^{-\ld j} =C M_0^\ld M_j^{-\ld}$ for some small $\ld > 0$
(to be determined later.) Then,  from \eqref{FW1}, we have
\begin{align} \label{A:subadd}
\mu \big(\|\mathbb{P}_{>M_0} u(\omega) \|_{\ft{b}^s_{p, 1} (\T)} > \eps)
& \leq \sum_{j = 1}^\infty \mu \big( \| \{\jb{n}^{s}n^{-1}
g_n(\omega)  \}_{|n| \sim M_j} \|_{\ell_n^{p}}  > \s_j \eps \big).
\end{align}

\noindent By interpolation and \eqref{A:decay}, we have
\begin{align*}
\| \{ & \jb{n}^s n^{-1} g_n \}_{|n| \sim M_j} \|_{\ell_n^{p}} \sim
M_j^{s-1} \| \{ g_n \}_{|n| \sim M_j} \|_{\ell_n^{p}} \leq  M_j^{s-1}
\| \{ g_n \}_{|n| \sim M_j} \|_{\ell_n^{2}}^\frac{2}{p}
  \| \{ g_n \}_{|n| \sim M_j} \|_{\ell_n^{\infty}}^\frac{p-2}{p} \\
& \leq  M_j^{s-1} \| \{ g_n \}_{|n| \sim M} \|_{\ell_n^{2}}
\Bigg(\frac{ \| \{ g_n \}_{|n| \sim M_j} \|_{\ell_n^{\infty}} } {\| \{
g_n \}_{|n| \sim M_j} \|_{\ell_n^{2}}} \Bigg)^\frac{p-2}{p} \leq
M_j^{s-1 -\dl \frac{p-2}{p}} \| \{ g_n \}_{|n| \sim M_j}
\|_{\ell_n^{2}}.
\end{align*}

\noindent Thus, if we have $\|\{\jb{n}^s n^{-1}g_n \}_{|n| \sim M_j}
\|_{\ell_n^{p}}  > \s_j \eps$,
 then we have
 $\| \{ g_n \}_{|n| \sim M_j} \|_{\ell_n^{2}}
 \gtrsim R_j $
 where $R_j := \s_j \eps M_j^\al$ with $\al := -s+1+\dl \frac{p-2}{p} $.
With $p = 2 + \theta$, we have $\al =  \frac{-(s-1)p +  \dl
\theta}{2 +  \theta} > \frac{1}{2}$ by taking $\dl$ sufficiently
close to $\frac{1}{2}$ since $-(s-1)p > 1$. Then, by taking $\ld >
0$ sufficiently small, $R_j = \s_j \eps M_j^\al = C \eps M_0^\ld
M_j^{\al-\ld} \gtrsim C \eps M_0^{\ld} M_j^{\frac{1}{2}+} $. By a
direct computation in the polar coordinates, we have
\begin{align*}
\mu\big(  \| \{ g_n \}_{|n| \sim M_j} \|_{\ell_n^{2}}  \gtrsim R_j
\big)
 \sim \int_{B^c(0, R_j)} e^{-\frac{1}{2}|g_n|^2} \prod_{|n| \sim M_j} dg_n
 \lesssim \int_{R_j}^\infty e^{-\frac{1}{2}r^2}  r^{2 \cdot \# \{|n| \sim M_j\} -1} dr.
\end{align*}

\noindent Note that, in the inequality, we have dropped the implicit
constant $\s(S^{2 \cdot \# \{|n| \sim M_j\} -1})$, a surface measure
of the $2 \cdot \# \{|n| \sim M_j\} -1$ dimensional unit sphere,
since $\s(S^n) = 2 \pi^\frac{n}{2} / \Gamma (\frac{n}{2}) \lesssim
1$. By the change of variable $t = M_j^{-\frac{1}{2}}r$, we have
$r^{2 \cdot \# \{|n| \sim M_j\} -2} \lesssim r^{4M_j} \sim
M_j^{2M_j} t^{4M_j}.$ Since $t > M_j^{-\frac{1}{2}} R_j = C \eps
M_0^{\ld} M_j^{0+}$, we have $M_j^{2M_j} = e^{2M_j \ln M_j} <
e^{\frac{1}{8}M_jt^2}$ and $ t^{4M_j} < e^{\frac{1}{8}M_jt^2}$ for
$M_0$ sufficiently large. Thus, we have $r^{2 \cdot \# \{|n| \sim
M_j\} -2} < e^{\frac{1}{4}M_jt^2} = e^{\frac{1}{4}r^2}$ for $ r >
R.$ Hence,  we have
\begin{align} \label{A:highfreq1}
\mu\big(   \| \{ g_n \}_{|n| \sim M_j} \|_{\ell_n^{2}}  \gtrsim R_j
\big) \leq C \int_{R_j}^\infty e^{-\frac{1}{4}r^2} r dr \leq
e^{-cR_j^2} = e^{-cC^2 M_0^{2\ld} M_j^{1+} \eps^2}.
\end{align}

\noindent From \eqref{A:subadd} and \eqref{A:highfreq1}, we have
\begin{align*}
\mu \big(\|\mathbb{P}_{>M_0}u\|_{\ft{b}^s_{p, 1} (\T)} > \eps) \leq
\sum_{j =1}^\infty e^{-cC^2 M_0^{1+2\ld+} (2^{j})^{1+} \eps^2} \leq
\tfrac{1}{2} \eps
\end{align*}

\noi by choosing $M_0$ sufficiently large as long as $(s-1)p<-1$.

When $p = \infty$, we have $s <1$. By repeating the computation with
\eqref{A:decay}, we see that if we have $\|\{\jb{n}^s n^{-1}g_n
\}_{|n| \sim M_j} \|_{\ell_n^{\infty}}  > \s_j \eps$, then we have
 $\| \{ g_n \}_{|n| \sim M_j} \|_{\ell_n^{2}}
 \gtrsim R_j $
 where $R_j := \s_j \eps M_j^\al$ with $\al := -s+1+\dl $.
Since $ -s+1 > 0$, we have $\al > \frac{1}{2}$ by taking $\dl$
sufficiently close to $\frac{1}{2}$. The rest follows exactly as
before.
\end{proof}

\subsection{Large deviation estimates for $\ft{b}^s_{p, \infty} (\T)$
at the endpoint $(s-1)p = -1$} \label{SUBSEC:Large}

Now, we show that the condition \eqref{msble} actually fails for
$\ft{b}^{s}_{p, \infty} (\T)$ for the endpoint case $(s-1) p = -1$; see
Remark \ref{Rmk5}. By the strong law of large numbers,
$X_j^{(p)}$ defined in \eqref{XJ} converges a.s.~to $c_p >0$. Then,
by Egoroff's theorem, there exists $E$ with $Pr(E) \geq \frac{1}{2}$
such that $X_j$ converges uniformly to $c_p$ on $E$. Thus, given
$\delta > 0$, there exists $J_0 \in \mathbb{N}$ such that
$Pr(\{\omega: \sup_{j \geq J} X_j^{(p)}(\omega) > c_p - \dl) \geq
\frac{1}{2}$ for any $J \geq J_0$. In view of \eqref{CON2}, this
shows that the condition \eqref{msble} does not hold once $\eps <
c_p$. In particular, Proposition \ref{PROP:FER} does not hold
automatically.

The remainder of this subsection is dedicated to the proof of
Theorem \ref{THM:end2} via a direct approach that bypasses the
assumption of abstract Wiener space. Specifically, we establish
that, for some $c = c(p)$ and all sufficiently large $K \geq K_p$,
the large deviation estimate
\begin{equation}\label{largedevi1}
Pr( \| u\|_{\ft{b}^s_{p, \infty} (\T)} > K) < e^{-cK^2},
\end{equation}

\noi also holds in the endpoint case $(s-1) p = -1$.

Let us first consider the case $p \leq 2$. By H\"older inequality
with $\frac{1}{r} = \frac{1}{p} - \frac{1}{2}$, we have
\begin{equation} \label{p20}
 \big\| 2^{-\frac{j}{p}} |g_n| \big\|_{\ell^p_{|n|\sim 2^j}}
\leq \| 2^{-\frac{j}{r}}\|_{\ell^r_{|n|\sim 2^j}}
\big\|2^{-\frac{j}{2}} |g_n| \big\|_{\ell^2_{|n|\sim 2^j}} \sim
\big\|2^{-\frac{j}{2}} |g_n| \big\|_{\ell^2_{|n|\sim 2^j}}.
\end{equation}

\noi i.e., we have $\| u\|_{\ft{b}^s_{p, \infty} (\T)} \lesssim \|
u\|_{\ft{b}^\frac{1}{2}_{2, \infty} (\T)}$. Hence, we only need to
consider the case $p = 2$ and $ s = \frac{1}{2}$. By definition of
the norm, we have
\begin{equation} \label{p21}
Pr(\| u\|_{\ft{b}^\frac{1}{2}_{2, \infty} (\T)} >K) \leq \sum_{j =
0}^\infty Pr\big( 2^{-j} \sum_{|n|\sim 2^j} |g_n|^2 > K^2\big).
\end{equation}

Let us now recall the so called {\it Cram\'er condition}: a sequence
$\{\xi_n\}$ of independent identically distributed (i.i.d.) random
variables is said to satisfy Cram\'er's condition if there exists
$\ld > 0$ such that
$$\varphi(\ld) = \mathbb{E}\big[
e^{\ld|\xi_1|}\big] < \infty.$$
If the condition holds, then we can
define {\it the Cram\'er transform}
$$H(a) = \sup_{\ld>0} \{ a \ld - \psi(\ld) \},$$
with $\psi(\ld) = \ln \varphi(\ld)$, and, for $a >
\mathbb{E}\big[\xi_1\big]$, we have
\begin{equation} \label{p22}
 Pr\Big(\frac{1}{N}\sum_{n=1}^N \xi_n \geq a \Big) \leq e^{-n H(a)};
\end{equation}

\noi see Shiryaev's book \cite{SHIR}. With $\xi_n = |g_n|^p$ (note
that $g_n$ is complex-valued), we see that Cram\'er's condition is
satisfied for $p=2$. Indeed, when $p = 2$, we have
\begin{align*}
\mathbb{E}\big[ e^{\ld |g_n|^2}\big] = \frac{1}{2\pi} \int e^{(\ld -
\frac{1}{2})|g_n|^2} d g_n = \frac{1}{2\pi(1-2\ld)}
\end{align*}

\noi for $\ld < \frac{1}{2}$. Then, $H(a) = \sup_{\ld
> 0} \{ a \ld + \ln (1-2\ld) + \ln (2\pi)\}$ has the maximum value
$\frac{a-2}{2} + \ln \frac{2}{a} + \ln (2\pi)$ at $\ld =
\frac{a-2}{2a}$. Then, from \eqref{p22}, we have
\begin{align} \label{B5}
Pr\big( 2^{-j} \sum_{|n|\sim 2^j} |g_n|^2 > K^2\big) < e^{-c 2^j
K^2}.
\end{align}

\noi Hence, we have $Pr(\| u\|_{\ft{b}^\frac{1}{2}_{2, \infty} (\T)} >K)
< e^{-c  K^2}$ in view of \eqref{p21}. This proves \eqref{largedevi1}
for $p\leq 2$.

Note that the Cram\'er condition no longer holds for $p > 2$. Thus,
we need another approach. There are known large deviation results
even when the Cram\'er's condition fails; see, for example, Saulis
and Nakas' work \cite{SN}. However, they do not seem to be directly
applicable to obtain \eqref{largedevi1}. Instead, we use the
hypercontractivity of the Ornstein-Uhlenbeck semigroup related to
products of Gaussian random variables. For the following discussion,
see the works of Kuo \cite{KUO2}, Ledoux-Talagrand \cite{Ledoux},
and Janson \cite{Janson}. A nice summary is given by Tzvetkov in
\cite[Section 3]{TZ3}.

In our discussion, we will use the Hermite polynomials $H_n(x)$. They
are defined by
\[ e^{tx -\frac{1}{2}t^2} = \sum_{n = 0}^\infty \frac{H_n(x)}{n!} t^n.\]

\noi The first three Hermite polynomials are: $H_0(x) = 1$, $H_1(x)
= x$, and $H_2(x) = x^2-1$.

Now, consider the Hilbert space $H = L^2(\R^m, \mu_m)$ with $d\mu_m
= (2\pi)^{-\frac{m}{2}} \exp(-{|x|^2}/{2})dx$, $x = (x_1, \dots,
x_m)\in \R^m$. We define a {\it homogeneous Wiener chaos of order
$n$} to be an element of the form $\prod_{j = 1}^m H_{n_j}(x_j)$, $n
= n_1 + \cdots + n_m$. Consider the Hartree-Fock operator $L = \Dl -
x \cdot \nabla$, which is the generator for the Ornstein-Uhlenbeck
semigroup. Then, by the hypercontractivity of the Ornstein-Uhlenbeck
semigroup $S(t) = e^{Lt}$, we have the following
\begin{lemma} \label{LEM:hyp}
Fix $q \geq 2$. Then, for every $u \in H$ and $t \geq
\frac{1}{2}\log(q-1)$, we have
\begin{equation}\label{hyp1}
\|S(t) u \|_{L^q(\R^m, \mu_m)}\leq \|u\|_{L^2(\R^m, \mu_m)}.
\end{equation}
\end{lemma}

\noi Note that \eqref{hyp1} holds, independent of the dimension $m$.
It is known that the eigenfunction of $L$ with eigenvalue $-n$
is precisely the homogeneous Wiener chaos of order $n$. Thus, we
have
\begin{lemma} \label{PROP:hyp}
Let $F(x)$ be a linear combination of homogeneous chaoses of order
$n$. Then, for $q \geq 2$, we have
\begin{equation} \label{hyp2}
\| F(x) \|_{L^q(\R^m, \mu_m)}\leq (q-1)^\frac{n}{2}
\|F(x)\|_{L^2(\R^m, \mu_m)}.
\end{equation}
\end{lemma}

\noi The proof is basically the same as in \cite[Propositions
3.3--3.5]{TZ3}. We only have to note that $F(x)$ is an eigenfunction
of $S(t) = e^{Lt}$ with eigenvalue $e^{-nt}$. Then, \eqref{hyp2}
follows from \eqref{hyp1} by evaluating \eqref{hyp1} at time $t =
\frac{1}{2} \log (q - 1)$.

Denote now by $\mathcal{K}_n $ the collection of the homogeneous
chaoses of order $n$. Given a homogeneous polynomial $P_n(x) =
P_n(x_1, \dots, x_m)$ of degree $n$, we define {\it the Wick ordered
monomial $:\!P_n(x)\!\!:$} to be its projection onto
$\mathcal{K}_n$. In particular, we have $:x_j^n: = H_{n}(x_j)$ and
$:\prod_{j = 1}^m x_j^{n_j}: = \prod_{j = 1}^m H_{n_j}(x_j)$ with $n
= n_1 + \cdots + n_m$.

Since the Fourier coefficients of Brownian motion involve complex
Gaussian random variables, let us consider the Wick ordering on them
as well. Let $g$ denote a standard complex-valued Gaussian random
variable. Then, $g$ can be written as $g = x+ iy$, where $x$ and $y$
are independent standard real-valued Gaussian random variables. Note
that the variance of $g$ is $\text{Var}(g)  = 2$. Next, we
investigate the Wick ordering on $|g|^{2n}$ for $n\in \mathbb{N}$,
that is, the projection of $|g|^{2n}$ onto $\mathcal{K}_{2n}$.

When $n = 1$, $|g|^2 = x^2 + y^2$ is Wick-ordered into
\[ :|g|^2: = (x^2 - 1) + (y^2-1) = |g|^2 - \text{Var}(g).\]

When $n = 2$, $|g|^4 = (x^2+ y^2)^2 = x^4 + 2 x^2 y^2 + y^4$ is
Wick-ordered into
\begin{align*}
:|g|^4: & = (x^4 -6 x^2 + 3) + 2(x^2 - 1)(y^2 - 1)
+ (y^4 -6 y^2 + 3)\\
& = x^4 + 2 x^2 y^2 + y^4
- 8 (x^2 + y^2) + 8 \\
& = |g|^4 - 4 \text{Var}(g) |g|^2 + 2 \text{Var}(g)^2,
\end{align*}

\noi where we used $H_4(x) = x^4 - 6x^2 + 3$.

When $n = 3$, $|g|^6 = (x^2+ y^2)^3 = x^6 + 3 x^4 y^2 +3x^2 y^4 +
y^6$ is Wick-ordered into
\begin{align*}
:|g|^6: & =  (x^6 - 15 x^4 + 45 x^2 - 15)
+ 3(x^4 - 6x^2 + 3)(y^2-1) \\
& \hphantom{X} + 3(x^2-1)(y^4 - 6y^2 + 3)
+ (y^6 - 15 y^4 + 45 y^2 - 15)\\
& = |g|^6 - 9\text{Var}(g)|g|^4 + 18 \text{Var}(g)^2|g|^2 - 6
\text{Var}(g)^3,
\end{align*}

\noi where we used $H_6(x) = x^6 - 15 x^4 + 45 x^2 - 15$.

In general, we have $:|g|^{2n}\!: \, \in \mathcal{K}_{2n}$.
Moreover, we have
\begin{equation}\label{Wick1}
:|g|^{2n}\!: \ = |g|^{2n} + \sum_{j = 0}^{n-1} a_j |g|^{2j}
 = |g|^{2n} + \sum_{j = 0}^{n-1} b_j :|g|^{2j}:  .
\end{equation}

\noi This follows from the fact that $|g|^{2n}$, as a polynomial in
$x$ and $y$ only with even powers, is orthogonal to any homogeneous
chaos of odd order, and it is radial, i.e., it depends only on $|g|^2
= x^2 + y^2$. Note that $:|g|^{2n}\!:$ can also be obtained from the
Gram-Schmidt process applied to $|g|^{2k}$, $ k= 0, \dots, n$ with
$\mu_2 = (2\pi)^{-1} \exp (-(x^2 + y^2)/2) dx dy$.

With these preliminaries, we are ready to return to the proof of the
large deviation estimate \eqref{largedevi1} for $p > 2$. Given $p >
2$, choose $k$ such that $p \leq 2k$. As in \eqref{p20}, by H\"older
inequality with $\frac{1}{r} = \frac{1}{p} - \frac{1}{2k}$, we have
\begin{equation}
 \big\| 2^{-\frac{j}{p}} |g_n| \big\|_{\ell^p_{|n|\sim 2^j}}
\lesssim \big\|2^{-\frac{j}{2k}} |g_n| \big\|_{\ell^{2k}_{|n|\sim
2^j}},
\end{equation}

\noi i.e., we have $\| u\|_{\ft{b}^s_{p, \infty} (\T)} \lesssim \|
u\|_{\ft{b}^{1-\frac{1}{2k}}_{{2k}, \infty} (\T)}$ for $(s-1)p = -1$.
Hence, it suffices to prove  \eqref{largedevi1} for $p = 2k$ and $s
= 1-\frac{1}{2k}$. Let
\begin{equation} \label{FJ}
F_j (\omega) = 2^{-j} \sum_{|n|\sim 2^j} |g_n(\omega)|^{2k}.
\end{equation}

\noi
 Then, we have
\[Pr (\sup_j |F_j| > K^{2k}) \leq \sum_{j = 0 }^\infty Pr( |F_j| > K^{2k}).\]

\noi Hence, it suffices to prove
\begin{equation} \label{largedevi2}
\sum_{j = 0 }^\infty Pr ( |F_j| > K^{2k}) < e^{-cK^2}.
\end{equation}

\noi By \eqref{Wick1}, write $F_j$ as a linear combination of
homogeneous chaoses of order $2\ell$, $\ell = 0, 1,\dots, k$, i.e., we
have $F_j = \sum_{\ell = 0}^k F_j^{(\ell)}$, where $F_j^{(\ell)}$ is the
component of $F_j$ projected onto $\mathcal{K}_{2\ell}$. Then, it
suffices to prove
\begin{equation} \label{largedevi3}
\sum_{j = 0 }^\infty Pr ( |F_j^{(\ell)}| > \tfrac{1}{k+1}K^{2k}) <
e^{-cK^2}
\end{equation}

\noi for each $\ell = 0, 1, \dots, k$. By choosing $K$ sufficiently
large, we see that \eqref{largedevi3} trivially holds for $\ell = 0$,
since $F_j^{(0)}$ is a constant independent of $j$ and thus the
left-hand side of \eqref{largedevi3} is 0 for large $K$. For $\ell \geq
1$, it follows from  Lemma \ref{PROP:hyp} that, for $q\geq 2$,
\begin{equation} \label{Wick2}
\| F_j^{(\ell)} \|_{L^q(\Omega)} \leq C_\ell q^\ell \|F_j^{(\ell)}
\|_{L^2(\Omega)} = C_\ell' 2^{-\frac{j}{2}} q^{\ell}
\end{equation}
where the constants $C_\ell$ and $C'_\ell$ are independent of $j$.

Let us now recall the following
\begin{lemma}[Lemma 4.5 in \cite{TZ3}] \label{LEM:hyp2}
Suppose that we have, for all $q \geq 2$,
\[ \|F(\omega)\|_{L^q(\Omega)} \leq C N^{-\al} q^{\frac{n}{2}}\]
for some $\al, N, C > 0$ and $n \in \mathbb{N}$. Then, there exist
$c$ and $C'$ depending on $C$ and $n$ but {\rm independent} of $\al$
and $N$ such that
\[Pr ( |F(\omega)| > \ld ) \leq C' e^{-c N^\frac{2\al}{n}\ld^{\frac{2}{n}}}.\]
\end{lemma}

Thus, from \eqref{Wick2} and Lemma \ref{LEM:hyp2} with $n = 2\ell$, $N
= 2^j$, $\al = \frac{1}{2}$, and $\ld = \frac{1}{k+1}K^{2k}$, we have
\begin{equation*}
Pr ( |F_j^{(\ell)}| > \tfrac{1}{k+1}K^{2k}) < e^{-c_k 2^{\frac{j}{2\ell}}
K^{\frac{2k}{\ell}}} < e^{-c_k 2^{\frac{j}{2\ell}} K^2}.
\end{equation*}
This establishes \eqref{largedevi3}, and hence \eqref{largedevi2}
and \eqref{largedevi1}.

\begin{remark} \rm
With $(s-1)p = -1$, we have
\[\mathbb{E}[\|u_j\|_{\ft{b}^s_{p, \infty} (\T)}^p] \sim \mathbb{E} [X^{(p)}_j] = c_p,\]
where $u_j = \mathbb{P}_{|n|\sim 2^j} u$, and $X^{(p)}_j$ is defined
in \eqref{XJ}. Also, note that
\[\mathbb{E} [F_j] = \mathbb{E}[F_j^{(0)}] = F_j^{(0)}= c_p.\]

Hence, it follows from the above computation for $\ell = 1, \dots, k$
that
\begin{align}\label{LG1}
 Pr \big( \big|\|u\|_{\ft{b}^s_{p, \infty} (\T)} -
c_p^{\frac{1}{p}}\big| > K\big)
& < e^{-cK^2} \\
\label{LG2} Pr \big( \big|\|\mathbb{P}_{|n| \geq 2^N}
u\|_{\ft{b}^s_{p, \infty} (\T)} - c_p^{\frac{1}{p}}\big| > K\big) & <
e^{-c2^\frac{N}{p}K^2}.
\end{align}

\noi In probability theory, large deviation estimates are commonly
stated as in the estimates \eqref{LG1} and \eqref{LG2}. However, in
applications to partial differential equations, it is more common to
encounter these estimate in the form \eqref{Large1}; see \cite{BO4,
BO5, BO7, BO8,BT1, BT3, nors, OH3, OH4, OHSBO}.

\end{remark}

\subsection{Large deviation estimates for $B_{p, \infty}^\frac{1}{2}$}

Lastly, we briefly discuss the large deviation estimates on the
Besov spaces $B^s_{p, q}$ with the endpoint regularity $s =
\frac{1}{2}$, $p < \infty$, and $ q= \infty$:
\begin{equation}\label{largedevi4}
Pr (\| u\|_{B^\frac{1}{2}_{p, \infty}} > K) < e^{-cK^2}
\end{equation}

\noi for some $c = c(p)$ and all sufficiently large $K \geq K_p$.
For the non-endpoint result, the reader is referred to \cite{ROY}.

For $p \leq 2$, \eqref{largedevi4} follows from \eqref{largedevi1}
once we note that $\| u\|_{B^\frac{1}{2}_{p, \infty}} \leq \|
u\|_{B^\frac{1}{2}_{2, \infty}} = \| u\|_{\ft{b}^\frac{1}{2}_{2,
\infty}}$. When $p > 2$, \eqref{largedevi2} does not follow from
\eqref{largedevi1} anymore. But, as in the proof of
\eqref{largedevi1}, it suffices to consider only the case $p = 2k$,
$k \geq 2$. The proof for a general even index $p$ involves lots of
unwieldy technicalities. In the following we will sketch the
argument for $p = 4$.

When $p = 4$, we have
\[ \| u\|_{B^\frac{1}{2}_{4, \infty}}^4
 = \sup_j \Big\|2^{-\frac{j}{2}} \sum_{|n|\sim 2^j} g_n
e^{int}\Big\|_{L^4_t}^4 = \sup_j \big( \I^{(2)}_j + \II^{(2)}_j +
\III^{(2)}_j\big), \]

\noi where $\I^{(2)}_j$, $ \II^{(2)}_j$, and  $\III^{(2)}_j$ are
defined in \eqref{B4}. In the following, we treat them separately.

First, note that
\begin{align*}
 \big\{ \omega: \, \I^{(2)}_j(\omega) >  2 K^4 \big\}
\subset \bigcup_{\ell = 1}^2 \bigg\{\omega: \, 2^{-j} \sum_{|n_\ell| \sim
2^j} |g_{n_\ell}(\omega)|^2 > K^2 \bigg\}.
\end{align*}

\noi Then, from \eqref{B5}, we have
 \begin{equation} \label{J1}
 Pr \big( \sup_j \I^{(2)}_j(\omega) >K^4\big) < e^{-cK^2}.
 \end{equation}

Next, note that $\III^{(2)}_j(\omega) >  K^4$ if and only if $F_j
(\omega) > 2^j K^4$, where $F_j(\omega)$ is defined in \eqref{FJ}
with $k = 2$. Hence, from \eqref{largedevi2}, we have
\begin{equation} \label{J2}
 Pr \big( \sup_j \III^{(2)}_j(\omega) >K^4\big) < e^{-cK^2}.
 \end{equation}

Lastly, by expanding the complex Gaussians into their real and
imaginary parts, it is not difficult to see that $\II^{(2)}_j$ is a
homogeneous Wiener chaos of order 4
since each term in the sum is a product of four independent real-valued
Gaussian random variables.
Then, it follows from  Lemma
\ref{PROP:hyp} that, for $q\geq 2$,
\begin{equation*}
\| \II^{(2)}_j \|_{L^q(\Omega)} \leq C q^2
\|\II^{(2)}_j\|_{L^2(\Omega)} = C' 2^{-\frac{j}{2}} q^{2}
\end{equation*}
where the constants $C$ and $C'$ are independent of $j$. Thus, from
Lemma \ref{LEM:hyp2}, we have $Pr ( |\II^{(2)}_j| > K^{4}) < e^{-c
2^{\frac{j}{4}}K^2}.$ This immediately implies
\begin{equation} \label{J3}
 Pr \big( \sup_j \big|\II^{(2)}_j(\omega)\big| >K^4\big) < e^{-cK^2}.
\end{equation}

\noi The large deviation estimate \eqref{largedevi4} follows from
\eqref{J1}, \eqref{J2}, and \eqref{J3}.

For a general even index $p$,
 one needs to repeat the above argument,
using \eqref{BP}. The estimates on $\I^{(k)}_j$ and $\II^{(k)}_j$
follow easily as before. In particular, note that $\II^{(k)}_j$ is a
homogeneous Wiener chaos of order $2k$. One can then estimate the
error terms by a combination of the arguments presented above.
However, the actual computation becomes lengthy, and thus we omit
details.

\appendix

\section{Brownian motion and Fourier-Wiener series}\label{ApA}
We present here a proof of the Fourier-Wiener series representation \eqref{2FW} of
the mean zero Brownian loop $u(t)$ on $[0, 2\pi)$.

Let $b(t)$ be the complex-valued Brownian motion on $\R_+$
and $\beta(t) = b(t) - t b(2\pi) /2\pi$ be the corresponding periodic Brownian loop on $[0, 2\pi)$.
Then, for $t \in [0, 2\pi)$, we have
\begin{align} \label{BMZ1}
d \beta = db - \frac{b(2\pi)}{2\pi} dt.
\end{align}

\noi
For $t \in [0, 2\pi)$,  we write them in terms of the Wiener integrals:
\[b(t) = \int_{\R_+} \chi_{[0, t)} d b(t')
\quad \text{and} \quad\beta(t) = \int_0^{2\pi} \chi_{[0, t)} d b(t')
- \frac{b(2\pi)}{2\pi}\int_0^{2\pi}\chi_{[0, t)} dt'  .\]

\noi
Given a (deterministic) periodic function $f$ on $[0, 2\pi)$, we write
$f = f_0 + f_1$ where $f_0 = \frac{1}{2\pi} \int_0^{2\pi} f(t) dt$
and $f_1 = f - f_0$.  Note that $\int_0^{2\pi} f_1(t) dt = 0$.

Now, define a conjugate\footnote{Instead of a linear operator,
we have a conjugate linear operator due to the complex-valued setting.} linear operator $T: L^2([0, 2\pi]) \to L^2 (\Omega)$ given by
\[T(f) (\omega) := \int_0^{2\pi} \cj{f} (t) d\beta(t; \omega).\]
Then,
from \eqref{BMZ1}, we have
\begin{align*}
T(f) & = \int_0^{2\pi} (\cj{f}_0 + \cj{f}_1)(t)  d b (t)
- \frac{b(2\pi)}{2\pi}\int_0^{2\pi}(\cj{f}_0 + \cj{f}_1)(t)  dt  \\
& = \int_0^{2\pi} \cj{f}_1 (t) d b (t)
+ \cj{f}_0 \int_0^{2\pi}   d b (t)
- b(2\pi)\cj{f}_0
= \int_0^{2\pi} \cj{f}_1 (t)  d b (t),
\end{align*}

\noi
where the last equality holds almost surely in $\omega$
since $b(0) = 0$ a.s.~Thus, we see that $T(f)$ is given by the Wiener integral  of the mean zero part $\cj{f}_1$.
Hence,
$T(f)$ is
a Gaussian random variable with mean 0 and variance $2 \|f_1\|^2_{L^2([0, 2\pi])}$.\footnote{The factor 2 appears
since we are considering the complex-valued Brownian motion.
See Kuo's book \cite{KUO2} for properties of the (real-valued) Wiener integrals.}
Namely, we have
\[ \|T(f)\|_{L^2(\Omega)} = \sqrt{2}\|f_1\|_{L^2([0, 2\pi])}.\]
Moreover, $T$ behaves like a (conjugate) unitary operator from $L^2([0, 2\pi])$ onto $L^2(\Omega)$, that is
\begin{equation} \label{BMZ2}
 \mathbb{E} \big[T(f)\cj{T(g)}\big] = 2\cj{\jb{f_1, g_1}}_{L^2}.
\end{equation}

On the one hand,  we have $\beta(t) = T(\chi_{[0, t)})$. On the other hand, we have
\[\chi_{[0, \tau)}(t) = \sum_{n\in \Z} a_n \cj{e}_n(t),\]

\noi
where $e_n(t) = (2\pi)^{-\frac{1}{2}}e^{int}$ and $a_n = \frac{e^{in\tau} - 1}{\sqrt{2\pi} i n}$.
Thus, we have
\begin{align} \label{BMZ3}
\beta(t) = T(\chi_{[0, t)}) = \sum_{n \ne 0} a_n g_n  =
 \sum_{n \ne 0}\frac{g_n}{\sqrt{2\pi}i n} \,  e^{int} - \sum_{n \ne 0}\frac{g_n}{\sqrt{2\pi} i n},
\end{align}

\noi
where $g_n = T(\cj{e}_n)$, $n \ne 0$,  is a Gaussian random variable with mean 0 and variance
$2\|e_n\|^2_{L^2} = 2$. (Note that $T(e_0) = 0$.)
Moreover, they are mutually independent from \eqref{BMZ2} and orthogonality of $\{e_n\}$.
Hence, the original Brownian motion $b(t)$ can be represented as
\[b(t) =
g_0 t +
 \sum_{n \ne 0}\frac{g_n}{\sqrt{2\pi}i n} \, e^{int} - \sum_{n \ne 0}\frac{g_n}{\sqrt{2\pi} i n} ,\]

\noi
where $g_ 0 = b(2\pi) / 2\pi$ is a Gaussian random variable with mean 0 and variance 2.
Also, by writing $g_0 =\frac{1}{2\pi} \int_0^{2\pi} db(t)$ a.s.~and
$g_n = \frac{1}{\sqrt{2\pi}} \int_0^{2\pi} e^{-int} db(t)$ for $n\ne 0$,
it follows that $g_0 $ is independent from $\{g_n\}_{n\ne 0}$.

Lastly, by subtracting the (spatial) mean of $\beta(t)$ over $[0, 2\pi)$ from \eqref{BMZ3}, we obtain
the Fourier-Wiener series representation of the mean zero periodic Brownian loop $u(t)$:
\begin{align} \label{BMZ4}
u(t) =  \sum_{n \ne 0}\frac{g_n}{\sqrt{2\pi}i n} \,  e^{int},
\end{align}

\noi
which is \eqref{2FW} up to the constant factor $\sqrt{2\pi}i $.

\begin{remark}\rm \label{RmkA}
Recall the following definition \cite{Janson}.
A Gaussian field on a Hilbert space $H$ is a (conjugate) linear isometry
(up to a multiplicative constant)
of $H$ into some Gaussian space.
In the previous discussion, we constructed such a
conjugate linear isometry in a concrete manner.

\end{remark}

\section{Brownian motion on the real line}\label{ApB}

In Remark \ref{Rmk1}, we defined the local-in-time versions of the time-frequency functions spaces in the following way.
Given an interval $I\subset \R$, we let $M^{p, q}_s (I)$  be the restriction of $M^{p, q}_s (\R)$ onto $I$
via
\begin{equation}\label{localized}
 \|u\|_{M^{p, q}_s(I)} = \inf \big\{ \|v\|_{M^{p, q}_s(\R)}: v = u \text{ on } I\big\}.
\end{equation}

\noi We define the local-in-time versions of other function spaces in an analogous manner.

In the following, we will show that, given a bounded interval $I$,  Theorem \ref{main} holds for the Brownian motion $b(t)$ on $\R$
in $M^{p, q}_s(I)$, $W^{p, q}_s(I)$, and $\ft{b}^s_{p, q}(I)$. Given the complex-valued Brownian motion $b(t)$ on $\R$,
it is known that $\sqrt{2a} b(t/a)$ and $b(t+t_0) - b(t_0)$,
$a\ne0$ and $t_0\in \R$,
are also Brownian motions.
Hence, it suffices to show that $b(t)$ is bounded or unbounded
a.s.~in
$M^{p, q}_s$, $W^{p, q}_s$, and $\ft{b}^s_{p, q}$ restricted to
some fixed bounded interval $I$.

We start by investigating the Brownian motion on local modulation spaces. We will repeatedly use the following proposition, the proof of which is
deferred to the end of this appendix.

\begin{proposition} \label{PROP:BMY1}
Let $f$ be a function on $\R$, and $\phi$ be a smooth cutoff function supported on $[0, 2\pi)$.
Then, for $1\leq p, q \leq \infty$ and $s \in \R $, we have

\noi
\textup{(a)}
\begin{equation} \label{BMY2}
\| \phi f  \|_{\mathcal{F}L^{s, q}(\T)} \lesssim \| f \|_{M^{p, q}_s(\R)}.
\end{equation}

\noi
\textup{(b)}
\begin{equation} \label{BMY2-2}
\| \phi f  \|_{M^{p, q}_s(\R)} \lesssim \| f \|_{\mathcal{F}L^{s, q}(\T)}.
\end{equation}

\end{proposition}

From \cite{RSTT}, we have
\begin{equation} \label{BMY1}
 \| f \|_{M^{p, q}_s(\R)} \sim \| f \|_{W^{p, q}_s(\R)}, \quad 1\leq p, q \leq \infty,
\end{equation}

\noi
for any function $f$ supported on a bounded interval $I$,
where the implicit constant depends on $p, q$, and $|I|$.

\subsection{Boundedness of Brownian motion}

Recall that we are free to select the fixed interval $I$. In this case, we take $I = [\frac{1}{2}\pi, \frac{3}{2}\pi)$.
Let $\phi$ be a smooth cutoff function supported on $[\frac{1}{4}\pi, \frac{7}{4}\pi)$
such that $\phi(t) = 1$ on $I$.
Also, let $\wt{\phi}$ be a smooth cutoff function supported on $[0, 2\pi)$
such that $\wt{\phi}(t) = 1$ on $\supp \phi$.
Note that $\phi = \phi \wt{\phi}$.
Then, by Proposition \ref{PROP:BMY1} (b), we have
\begin{align} \label{BMY4}
\|f\|_{M^{p, q}_s(I)} & \leq \|\phi f\|_{M^{p, q}_s(\R)}
= \|\phi \wt{\phi}f\|_{M^{p, q}_s(\R)}
\lesssim  \|\wt{\phi}f\|_{\mathcal{F}L^{s, q}(\T)} \notag \\
& \leq \|f\|_{\mathcal{F}L^{s, q}(\T)}
\sim \|f\|_{M^{p, q}_s(\T)},
\end{align}

\noi
where the first inequality follows from the definition of the localized space
and the last inequality follows from the boundedness of the multiplication
by a smooth function (supported on $[0, 2\pi)$.)
By \eqref{BMY1} and \eqref{BMY4}, we also have
\begin{align} \label{BMY5}
\|f\|_{W^{p, q}_s(I)} \lesssim \|f\|_{\mathcal{F}L^{s, q}(\T)}
\sim \|f\|_{W^{p, q}_s(\T)}.
\end{align}

Now,
let $\beta(t)$ be the periodic part of $b(t)$ on $I$,
i.e.
\[ \beta (t) = b(t) - b(\tfrac{1}{2}\pi) + \frac{b(\tfrac{3}{2}\pi) - b(\tfrac{1}{2}\pi)}{\pi}(t-\tfrac{1}{2}\pi).\]
Then, we have
\begin{equation} \label{BMY8}
\|b\|_{M^{p, q}_s(I)}
\leq \|\beta\|_{M^{p, q}_s(I)} + C_\omega(I),
\end{equation}

\noi
where $C_\omega(I) <\infty$ a.s.
Then, from Theorem \ref{main} with \eqref{BMY4} and \eqref{BMY8},
we have
\begin{align*}
\|b\|_{M^{p, q}_s(I)}
& \leq \|\beta\|_{M^{p, q}_s(I)} + C_\omega(I)
\lesssim  \|\phi \beta\|_{M^{p, q}_s(\R)} + C_\omega(I)\\
& \lesssim \|\beta\|_{M^{p, q}_s(\T)} + C_\omega(I) <\infty,\quad \text{a.s.}
\end{align*}

\noi
for $(s-1) q<-1$ with $q <\infty$ (and $s<1$ when $q = \infty$.)
The same boundedness result holds in $W^{p, q}_s(I)$
under the same condition.
%
%\noi
%From {\bf eqref}, \eqref{BMY1}, and \eqref{BMY2}, it follows that
%\begin{equation} \label{BMY3}
%\| u \|_{\mathcal{F}L^{s, q}(\R)} \sim \| u \|_{W^{p, q}_s(\R)} \sim \|u\|_{FL^{s, q}(\T)}
%\end{equation}
%
%\noi
%for a function $u$ supported on $[0, 2\pi)$.

\subsection{Unboundedness of Brownian motion}

In the following, we establish unboundedness of the Brownian motion on $\R$
in the function spaces restricted to bounded intervals on $\R$.
Due to the definition \eqref{localized} of the function spaces restricted to
an interval, it is more difficult to establish a lower bound on the localized norms.
(In order to to establish an upper bound on the localized norms, it suffices to
show an upper bound on a single representation - see \eqref{BMY4},
whereas we need to show that a uniform lower bound exists for {\it all}
representations  in establishing a lower bound on the localized norms.)
The trick we use is the following.
Given a periodic function of period $T$, we establish a lower bound on the norms restricted to the interval
of length $2T$.

Take $I = [0, \pi)$ such that $2I = [0, 2\pi) = \T$.
Let $\phi_1$ be a smooth cutoff function supported on $[\frac{1}{2}\pi - \eps, \frac{3}{2}\pi+\eps)$
with $\phi(t) = 1$ on $[\frac{1}{2}\pi + \eps, \frac{3}{2}\pi-\eps)$
for some small $\eps > 0$
such that $\phi_1 + \phi_2 \equiv 1$ on $\T$,
where $\phi_2(t) := \phi_1(t- \pi)$ on the periodic domain $\T$.
Also, let $\wt{\phi}_1$ be a smooth cutoff function supported  on $\T$
such that $\wt{\phi}_1(t) = 1$ on $\supp \phi_1$.
Note that $\phi_1 = \phi_1 \wt{\phi}_1$.

Let $f$ be a periodic function of period $\pi$.
Then, by Proposition \ref{PROP:BMY1} (a), we have
\begin{align*}
\|f\|_{M^{p, q}_s(\T)} & \sim \|f\|_{\mathcal{F}L^{s, q}(\T)}
= \|(\phi_1 + \phi_2) f\|_{\mathcal{F}L^{s, q}(\T)}
\leq \|\phi_1  f\|_{\mathcal{F}L^{s, q}(\T)}+\|\phi_2 f\|_{\mathcal{F}L^{s, q}(\T)}\\
\intertext{which is, noting that $\phi_2 (t) f(t) = \phi_1(t-\pi) f(t-\pi)$,}
& =2 \|\phi_1  \wt{\phi}_1 f\|_{\mathcal{F}L^{s, q}(\T)}
\lesssim \|\wt{\phi}_1 f\|_{M^{p, q}_s(\R)}.
\end{align*}

\noi
Let $g$ any function on $\R$ such that $g = f$ on $2I = [0, 2\pi)$.
Then, we have $\wt{\phi}_1 f = \wt{\phi}_1 g$.
By the boundedness of the multiplication by a smooth function, we obtain
\begin{align*}
\|f\|_{M^{p, q}_s(\T)}
\lesssim \|\wt{\phi}_1 g\|_{M^{p, q}_s(\R)}
\lesssim \| g\|_{M^{p, q}_s(\R)},
\end{align*}

\noi
for any extension $g$.
Hence, we have
\begin{align}
\|f\|_{M^{p, q}_s(\T)}
\lesssim \| f\|_{M^{p, q}_s(2I)},
\end{align}

\noi
for a periodic function $f$ of period $\pi$.
By the translation invariance property of the norm, we obtain
\begin{align} \label{BMY6}
\|f\|_{M^{p, q}_s(\T)}
\lesssim \| f\|_{M^{p, q}_s(2I)}
\leq \| f\|_{M^{p, q}_s(I)}  +\| f\|_{M^{p, q}_s(I+\pi)}
\lesssim  \| f\|_{M^{p, q}_s(I)}.
\end{align}

\noi
By \eqref{BMY1} and \eqref{BMY6}, we also have
\begin{align} \label{BMY7}
\|f\|_{W^{p, q}_s(\T)}
\lesssim \| f\|_{W^{p, q}_s(I)},
\end{align}

\noi
for a periodic function $f$ of period $\pi$.

Now,
let $\beta(t)$ be the periodic part of $b(t)$ on $[0, \pi)$,
i.e. $ \beta (t) = b(t)  - t b(\pi)/\pi.$
Then, we have
\begin{equation} \label{BMY9}
\|b\|_{M^{p, q}_s(I)}
\geq \|\beta\|_{M^{p, q}_s(I)} - C_\omega(I),
\end{equation}

\noi
where $C_\omega(I) <\infty$ a.s.
Then, from Theorem \ref{main} with \eqref{BMY6} and \eqref{BMY9},
we have
\begin{align*}
\|b\|_{M^{p, q}_s(I)}
 \geq \|\beta\|_{M^{p, q}_s(I)} - C_\omega(I)
\gtrsim  \| \beta\|_{M^{p, q}_s(\T)} - C_\omega(I)
 = \infty, \quad \text{a.s.}
\end{align*}

\noi
for $(s-1) q\geq -1$ with $q <\infty$ (and $s\geq 1$ when $q = \infty$.)
The same unboundedness result holds in $W^{p, q}_s(I)$
under the same condition.

\medskip

A similar argument with the following proposition shows the boundedness or unboundedness of the Brownian motion
in $\ft{b}^s_{p, q}(I)$ for any bounded interval $I$.

\begin{proposition} \label{PROP:BMA1}
 Let $f$ be a function on $\R$, and $\phi$ be a smooth cutoff function supported on $[0, 2\pi)$.
Then, for $1\leq p, q \leq \infty$ and $s \in \R $, we have

\noi
\textup{(a)}
\begin{equation} \label{BMA1}
\| \phi f  \|_{\ft{b}^s_{p, q}(\T)} \lesssim \| f \|_{\ft{b}^s_{p, q}(\R)}.
\end{equation}

\noi
\textup{(b)}
\begin{equation} \label{BMA2}
\| \phi f  \|_{\ft{b}^s_{p, q}(\R)} \lesssim \| f \|_{\ft{b}^s_{p, q}(\T)}.
\end{equation}
\end{proposition}

We conclude this appendix and our paper by presenting the proofs of
Propositions \ref{PROP:BMY1}
and  \ref{PROP:BMA1}.

\begin{proof}[Proof of Proposition \ref{PROP:BMY1}]
(a)
First, recall the support function $\psi$ in the definition of the modulation spaces,
i.e.~$\psi$ is a smooth cutoff function supported on $[-1, 1]$
such that $\sum_{n\in\Z} \psi(\xi - n) \equiv 1$.

Now, let $\eta(\xi)$ be a smooth cutoff function supported on $[-2, 2]$
such that $\eta(\xi) = 1$ on $[-1, 1]$.
Then, we have
$\psi(\xi) \eta(\xi) = \psi(\xi)$ for $\xi \in \R$.
For $n\in \Z$, we let $T_n$ and $M_n$ denote the translation and modulation operators respectively, that is,
$T_n\eta (\xi) = \eta(\xi - n)$ and $M_n\phi (y) = \phi(y)  e^{iny}$.
Then, we have
\begin{align*}
\ft{\phi f}(n) & = \int_0^{2\pi} \phi(y) f(y) e^{-iny} dy
= \sum_{m\in \Z}\int_\R \psi(D-m) f (y) \cj{M_n\phi (y)} dy\\
& = \sum_{m\in \Z}\int_\R \psi(\xi-m) \ft{f} (\xi) \cj{\ft{M_n \phi}(\xi)} d\xi
= \sum_{m\in \Z}\int_\R \psi(\xi-m) \ft{f} (\xi) \cj{T_m\eta (\xi)\ft{M_n\phi}(\xi)} d\xi\\
& = \sum_{m\in \Z}\int_\R \psi(D-m)f(y) \cj{\invft{T_m\eta}*M_n\phi (y)} dy.
\end{align*}

First, we prove \eqref{BMY2} for $s\geq 0$.
For $m\ne n$, repeated integration by parts gives
\begin{align*}
|\invft{T_m\eta}*M_n\phi (x)| & = \bigg|\int_\R \invft{\eta}(x - y) e^{im(x-y)} \phi(y)e^{iny} dy\bigg|
= \bigg|\int_\R \invft{\eta}(x - y)  \phi(y)e^{i(n-m)y} dy\bigg|\\
& \lesssim \frac{1}{\jb{n-m}^{s+2}}\bigg| \int_\R \frac{1}{\jb{x-y}^{2}} \wt{\phi}(y) dy \bigg|
\end{align*}

\noi
for some $\wt{\phi}\in\mathcal{S}(\R)$.
Thus, we have
\begin{equation} \label{BMX1}
\|\invft{T_m\eta}*M_n\phi \|_{L^{p'}(\R)} \lesssim \frac{1}{\jb{n-m}^{s+2}},\quad m \ne n,
\end{equation}

\noi
where $p'$ is the H\"older conjugate exponent of $p \in [1, \infty]$.
By H\"older inequality and \eqref{BMX1}, we have
\begin{align}
|\ft{\phi f}(n)|
& = \sum_{m\in \Z} \|\psi(D-m)f\|_{L^p(\R)} \|\invft{T_m\eta}*M_n\phi (y)\|_{L^{p'}(\R)} \notag \\
& \lesssim  \sum_{m\in \Z} \frac{1}{\jb{n-m}^{s+2}}\|\psi(D-m)f\|_{L^p(\R)}. \label{BMX2}
\end{align}

\noi
By $\jb{n}^s \lesssim \jb{n-m}^s \jb{m}^s$ for $s \geq 0$, \eqref{BMX2}
and Young's inequality, we obtain
\begin{align*}
\|\phi f\|_{\mathcal{F}L^{s, q}(\T)}
& = \|\jb{n}^s \ft{\phi f}(n)\|_{\ell^q(\Z)}
\lesssim \bigg\|\sum_{m\in \Z} \frac{1}{\jb{n-m}^{2}} \jb{m}^s \|\psi(D-m)f\|_{L^p(\R)}\bigg\|_{\ell^q_n(\Z)}\\
& \lesssim \big\| \jb{n}^s \|\psi(D-n)f\|_{L^p(\R)}\big\|_{\ell^q_n(\Z)}
\lesssim \|f\|_{M^{p, q}_s(\R)}.
\end{align*}

Next, assume $s< 0$.
In this case, we have $\jb{n}^s  \jb{n-m}^s \lesssim \jb{m}^s$.
By repeating the previous argument, we have
\begin{equation*}
\|\invft{T_m\eta}*M_n\phi\|_{L^{p'}(\R)} \lesssim \frac{1}{\jb{n-m}^{-s+2}},\quad m \ne n.
\end{equation*}

\noi
Then, the rest follows as before.

\medskip
\noi
(b) Recall that $\phi$ is a smooth cutoff function supported on $[0, 2\pi)$.
Let $n\in \Z$. Then, we have
\begin{align}
\psi(D-n) (\phi f)(x)
& = \int_\R \psi(\xi - n) \ft{\phi f}(\xi) e^{i\xi x} d\xi
= \int_\R \psi(\xi - n) \int_0^{2\pi}\phi(y) f(y) e^{-i \xi y}d y \, e^{i\xi x} d\xi\notag \\
& = \int_0^{2\pi} \phi(y) f(y) \int_\R \psi(\xi - n) e^{i\xi(x-y)} d\xi dy\notag \\
& = e^{inx} \int_0^{2\pi} f(y) \phi(y) \invft{\psi}(x-y) e^{-iny}dy \notag \\
& = e^{inx} \int_0^{2\pi} \sum_{m \in \Z} \ft{f}(n-m)e^{i(n-m)y} \phi(y) \invft{\psi}(x-y) e^{-iny}dy \notag \\
& = e^{inx} \sum_{m \in \Z} \ft{f}(n-m) F_m(x), \label{BMX3}
\end{align}

\noi
where $\ft{f}(n-m)$ are the Fourier coefficients of $f$ (as a function on $\T$)
at $n-m$ and
\[ F_m( x) = \int_\R \phi(y) \invft{\psi}(x-y) e^{-imy}dy.\]

First, we prove \eqref{BMY2-2} for $s\geq 0$.
By Young's inequality, we have
\begin{equation}\label{BMX4}
 \|F_0(x)\|_{L^p (\R)} = \bigg\|\int_\R \phi(y) \invft{\psi}(x-y)dy \bigg\|_{L^p_x(\R)}
\lesssim 1.
\end{equation}

\noi
For $m\ne 0$, by repeatedly integrating by parts as before, we obtain
\begin{align*}
|F_m(x)| &
 \lesssim \frac{1}{\jb{m}^{s+2}}\bigg| \int_\R \frac{1}{\jb{x-y}^{2}} \wt{\phi}(y) dy \bigg|
\end{align*}

\noi
for some $\wt{\phi}\in\mathcal{S}(\R)$.
Hence, by Young's inequality, we have
\begin{equation}\label{BMX5}
 \|F_m(x)\|_{L^p(\R)}  \lesssim \frac{1}{\jb{m}^{s+2}}.
\end{equation}

\noi
From \eqref{BMX3}, \eqref{BMX4}, and \eqref{BMX5}, we have
\begin{align}\label{BMX6}
\|\psi(D-n) f(x)\|_{L^p(\R)}
\lesssim \sum_{m \in \Z} |\ft{f}(n-m)| \frac{1}{\jb{m}^{s+2}}.
\end{align}

\noi
By $\jb{n}^s \lesssim \jb{n-m}^s \jb{m}^s$ for $s \geq 0$, \eqref{BMX6}
and Young's inequality,
we obtain
\begin{align*}
\|\phi f\|_{M^{p, q}_s(\R)} & = \big\|\jb{n}^s \|\psi(D-n) f(x)\|_{L^p(\R)}\big\|_{\ell^q_n(\Z)}
\lesssim \bigg\| \sum_{m \in \Z} \jb{n-m}^s |\ft{f}(n-m)| \frac{1}{\jb{m}^{2}} \bigg\|_{\ell^q_n(\Z)}\\
& \lesssim \|\jb{n}^s \ft{f}(n)\|_{\ell^q_n(\Z)} = \|f\|_{\mathcal{F}L^{s, q}(\T)}.
\end{align*}

\noi
The modification for $s<0$ follows as in Part (a).
\end{proof}

\begin{proof}[Proof of Proposition \ref{PROP:BMA1}]
(a) Let $\{\varphi_j\}_{j = 0}^\infty$ be the support function in the definition of
$\ft{b}^s_{p, q}(\R)$ (and the Besov space), i.e.
$\varphi_0, \varphi \in \mathcal{S}(\R)$
such that $\supp \varphi_0 \subset \{|\xi| \leq 2\}$
and $\supp \varphi \subset \{\frac{1}{2} \leq |\xi| \leq 2\}$
with $\varphi_0(\xi) + \sum_{j = 1}^\infty \varphi_j(\xi) \equiv 1$,
where $\varphi_j(\xi) = \varphi(2^{-j}\xi)$.

Let $\eta_0(\xi)$ be a smooth cutoff function supported on $\{ |\xi| \leq 3\}$
such that $\eta_0(\xi) = 1$ on $\{|\xi| \leq 2\}$.
For $j \geq 1$, let $\eta_j(\xi) = \eta(2^{-j}\xi)$,
where $\eta$ is a smooth cutoff function supported on $\{ \frac{1}{4} \leq |\xi| \leq 3\}$
such that $\eta(\xi) = 1$ on $\{\frac{1}{2}\leq |\xi| \leq 2\}$.
Then, we have $\varphi_j(\xi) \eta_j(\xi) = \varphi_j(\xi)$
for $\xi \in \R$ and $j\geq 0$.
Also, recall the notation $M_n\phi (y) = \phi(y) e^{iny}$.
Then, we have
\begin{align}
|\ft{\phi f}(n)| & = \bigg|\int_0^{2\pi} \phi(y) f(y) e^{-iny} dy\bigg|
= \bigg|\sum_{j = 0}^\infty\int_\R \phi(y) \varphi_j (D) f(y) e^{-iny} dy\bigg|\notag \\
& = \bigg|\sum_{j = 0}^\infty\int_\R  \varphi_j (\xi) \ft{f}(\xi) \cj{\eta_j (\xi) \ft{M_n\phi} (\xi)} dy\bigg|
\leq \sum_{j = 0}^\infty \|\varphi_j \ft{f}\|_{L^p_\xi(\R)}
\|\eta_j \ft{M_n\phi}\|_{L^{p'}_\xi(\R)}. \label{BMA5}
\end{align}

\noi
Since $\ft{\phi} \in \mathcal{S}(\R)$, we have
\begin{align*}
\|\eta_0 \ft{M_n\phi}\|_{L^{p'}_\xi(\R)}
& = \bigg(\int |\eta_0(\xi) \ft{\phi}(\xi - n)|^{p'} d\xi\bigg)^\frac{1}{p'}
\lesssim  c_{0, n},
\intertext{and}
\|\eta_j \ft{M_n\phi}\|_{L^{p'}_\xi(\R)}
& = \bigg(\int |\eta(2^{-j} \xi) \ft{\phi}(\xi - n)|^{p'} d\xi\bigg)^\frac{1}{p'}\\
& \lesssim  c_{j, n} \bigg(\int_{|\xi|\sim 2^j}  \frac{1}{\jb{\xi - n}^2} d\xi\bigg)^\frac{1}{p'}
\lesssim c_{j, n}, \quad j\geq 1,
\end{align*}

\noi
where
\begin{equation}c_{j, n} = \begin{cases} 1, & |j-k|\leq 2,\\
\jb{2^j - n}^{-2-s}, & |j-k|\geq 3 \text{ and } s \geq 0,\\
\jb{2^j - n}^{-2+s}  & |j-k|\geq 3 \text{ and } s < 0,
\end{cases}
\qquad \text{with } |n| \sim 2^k. \label{BMA3}
\end{equation}

From \eqref{BMA5}, we have
\begin{align} \label{BMA4}
\|\phi f \|_{\ft{b}^s_{p, q}(\T) }
= \big\| \| \jb{n}^s \ft{\phi f}(n)\|_{\ell^p_{|n|\sim 2^k} }\big\|_{\ell^q_k}
\lesssim  \bigg\| \Big\|2^{ks} \sum_{j = 0}^\infty c_{j, n}\|\varphi_j \ft{f}\|_{L^p_\xi(\R)} \Big\|_{\ell^p_{|n|\sim 2^k}}\bigg\|_{\ell^q_k}.
\end{align}

\noi
The contribution from  $|j - k|\leq 2$ is easily estimated by
\[ \lesssim \big\|2^{ks} \|\varphi_k \ft{f}\|_{L^p_\xi(\R)} \big\|_{\ell^q_k}
\sim \| f \|_{\ft{b}^s_{p, q}(\R) } .\]

When $j \geq k + 3$,
we have $c_{j, n} \sim 2^{-(2+s)j}$ for $|n| \sim 2^k$.
Then, the contribution to \eqref{BMA4} in this case is estimated by
\begin{align*}
\bigg\| \Big\|2^{ks} \sum_{j \geq k + 3 }&  c_{j, n}\|\varphi_j \ft{f}\|_{L^p_\xi(\R)}
\Big\|_{\ell^p_{|n|\sim 2^k}}\bigg\|_{\ell^q_k} \\
& \lesssim
\bigg\|
\big\|2^{js}\|\varphi_j \ft{f}\|_{L^p_\xi(\R)} \big\|_{\ell^q_j}
\Big\| \Big( \sum_{j \geq k + 3}2^{(k-j)sq'} 2^{-q'(2+s)j} \Big)^{\frac{1}{q'}}
\Big\|_{\ell^p_{|n|\sim 2^k}}\bigg\|_{\ell^q_k}\\
& \lesssim \| f \|_{\ft{b}^s_{p, q}(\R) }\,
\bigg(\sum_k \sum_{|n|\sim 2^k}
 \sum_{j \geq k + 3}2^{(k-j)s} 2^{-(2+s)j} \bigg)
\lesssim \| f \|_{\ft{b}^s_{p, q}(\R) },
\end{align*}

\noi
where we used H\"older inequality in $j$ in the first inequality.

When $j \leq k - 3$,
we have $c_{j, n} \sim 2^{-(2+s)k}$ for $|n| \sim 2^k$.
Then, the contribution to \eqref{BMA4} in this case is estimated by
\begin{align*}
\bigg\| \Big\|2^{ks} \sum_{j \leq k - 3 }&  c_{j, n}\|\varphi_j \ft{f}\|_{L^p_\xi(\R)}
\Big\|_{\ell^p_{|n|\sim 2^k}}\bigg\|_{\ell^q_k} \\
& \lesssim
\bigg\|
\big\|2^{js}\|\varphi_j \ft{f}\|_{L^p_\xi(\R)} \big\|_{\ell^q_j}
\Big\| \Big( \sum_{j \leq k - 3}2^{(k-j)sq'} 2^{-q'(2+s)k} \Big)^{\frac{1}{q'}}
\Big\|_{\ell^p_{|n|\sim 2^k}}\bigg\|_{\ell^q_k}\\
& \lesssim \| f \|_{\ft{b}^s_{p, q}(\R) }\,
\bigg(\sum_k \sum_{|n|\sim 2^k}
 \sum_{j \leq k - 3}2^{(k-j)s} 2^{-(2+s)k} \bigg)
\lesssim \| f \|_{\ft{b}^s_{p, q}(\R) }.
\end{align*}

\noi
Using \eqref{BMA3}, a similar computation yields \eqref{BMY1}
for $s< 0$.

\medskip

\noi
(b) Consider again $\phi$ to be a smooth cutoff function supported on $[0, 2\pi)$.
Then, we have
\begin{align}
\ft{\phi f}(\xi) & = \int_\R \phi(y) f(y) e^{-i\xi y} dy \notag \\
& = \int_0^{2\pi} \phi(y) f(y) e^{-iny} dy +
\int_0^{2\pi}  f(y)\phi_{n, \xi} (y)  e^{-iny} dy  \notag \\
& = \ft{\phi f}(n) +\ft{\phi_{n, \xi} f}(n), \label{BMA6}
\end{align}

\noi
where $\phi_{n, \xi} (y) = \phi(y) (e^{-i(\xi-n) y} - 1)$.

First, we prove \eqref{BMA2} for $s\geq 0$.
For $n \in \Z$, let $I_n = [n-1, n+1]$.
Then, for $\xi \in I_n$, we have
\begin{align}
|\ft{\phi_{n, \xi}}(m)| & = \bigg|\int_0^{2\pi} \phi(y) (e^{-i(\xi-n) y} - 1) e^{-imy} dy\bigg|
 = \bigg|\int_\R \phi(y) (e^{-i(\xi-n) y} - 1) e^{-imy} dy\bigg| \notag\\
& = |\ft{\phi}(\xi - n + m) - \ft{\phi}(m)|
\leq \sup_{\zeta \in [m-1, m+1]}|\partial_\zeta \ft{\phi}(\zeta)| |\xi-n|
\lesssim \jb{m}^{-2-s}. \label{BMA7}
\end{align}

\noi
Also, we have $|\ft{\phi}(m)| \lesssim \jb{m}^{-2-s}.$

From \eqref{BMA6} and \eqref{BMA7}, we have
\begin{align}
\|\phi f \|_{\ft{b}^s_{p, q}(\R) }
& = \big\| 2^{js}\|\varphi_j(\xi) \ft{\phi f}(\xi)\|_{L^p_\xi(\R)} \big\|_{\ell^q_j}
\lesssim \bigg\| 2^{js}  \bigg( \sum_{|n|\sim 2^j} \int_{\xi \in I_n}
|\ft{\phi f}(\xi)|^p d\xi\bigg)^\frac{1}{p}  \bigg\|_{\ell^q_j} \notag \\
& \lesssim
\bigg\| 2^{js}  \bigg( \sum_{|n|\sim 2^j}
\Big|\sum_m \ft{f}(n-m) \jb{m}^{-2-s}\Big|^p
\bigg)^\frac{1}{p}  \bigg\|_{\ell^q_j}. \label{BMA8}
\end{align}

\noi
Note that, for $|n| \sim 2^j$,  (i) $|n- m | \sim 2^j$ implies $|m| \lesssim 2^j$,
(ii) $|n-m| \ll 2^j$ implies $|m| \sim 2^j$, and
(iii) $|n-m| \gg 2^j$ implies $|m| \gg 2^j$.
First, we consider the case (i).
In this case, the contribution to \eqref{BMA8} is estimated by
\[ \lesssim
\bigg\| 2^{js}
\| \ft{f}(n)\|_{\ell^p_{|n|\sim2^j}}
\bigg( \sum_{|m| \lesssim  2^j} \jb{m}^{-2-s} \bigg)  \bigg\|_{\ell^q_j}
\lesssim \| f \|_{\ft{b}^s_{p, q}(\T) }.
\]

Next, we consider the case (ii).
By H\"older inequality, the contribution to \eqref{BMA8} is estimated by
\begin{align*} \lesssim
\bigg\| 2^{js}
\| \ft{f}(n)\|_{\ell^p_{|n|\ll2^j}}
\bigg( \sum_{|m| \sim  2^j} \jb{m}^{-2-s} \bigg)  \bigg\|_{\ell^q_j}
\lesssim
\bigg\|
j^\frac{1}{q'} 2^{-j}\,
\big\| \| \ft{f}(n)\|_{\ell^p_{|n|\sim 2^k}} \big\|_{\ell^q_k}
 \bigg\|_{\ell^q_j}\lesssim \| f \|_{\ft{b}^s_{p, q}(\T) }.
\end{align*}

Lastly, we consider the case (iii).
In this case, we have $|n-m|\sim |m|$.
By H\"older inequality, the contribution to \eqref{BMA8} is estimated by
\begin{align*}
& \lesssim
\bigg\| 2^{js}
\sum_{|m|\sim 2^k \gg 2^j} 2^{-(2+s)k} \| \ft{f}(n)\|_{\ell^p_{|n|\sim2^k}}
 \bigg)  \bigg\|_{\ell^q_j}\\
& \lesssim
\big\|2^{ks} \| \ft{f}(n)\|_{\ell^p_{|n|\sim 2^k}} \big\|_{\ell^q_k}
\bigg\|
\bigg(\sum_{k \geq j} 2^{-(2+s)q'k}\bigg)^\frac{1}{q'}
 \bigg\|_{\ell^q_j}\lesssim \| f \|_{\ft{b}^s_{p, q}(\T) }.
\end{align*}

\noi
When $s< 0$, we can instead choose
$|\ft{\phi_{n, \xi}}(m)|, |\ft{\phi}(m)| \lesssim \jb{m}^{-2+s}$  in \eqref{BMA7}.
The modification is straightforward, and thus we omit the details.
\end{proof}

\end{document}